\documentclass{amsart}
\newtheorem{thm}{Theorem}[section]
\newtheorem{prp}[thm]{Proposition}

\newtheorem{lma}[thm]{Lemma}
\theoremstyle{definition}
\newtheorem{dfn}[thm]{Definition}
\newtheorem{exa}{Example}

\def \a{a}
\def \b{b}
\def \ab{\bar\a}
\def \c{c}
\def \m{\pi}  
\def \ph{\varphi}
\def \sg{\sigma}

\def \A{\mbox{\bf A}}
\def \Ah{\A^\h}
\def \An{\Ah_n}
\def \Azero{\Ah_0}
\def \cycl{\operatorname {cycl}}
\def \Lb{[\cdot\,,\cdot]}
\def \Pb{\{\cdot\,,\cdot\}}
\def \SV{{\mathcal S}(V)}
\def \SVh{\SV^h}
\def \SVn{\SVh_n}
\def \s#1{\mbox{$(-1)^{#1}$}}
\def \sgn{\operatorname{sgn}}

\def \TV{{\mathcal T}(V)}
\def \TVh{\TV^h}
\def \TVn{\TVh_n}
\def \U{{\mathcal U}(V)}
\def \Uh{\U^h}
\def \Un{\Uh_n}
\def \Uhs{\Uh_\star}
\def \Jhs{\Jh_\star}
\def \Jh{J^\h}
\def \Jn{\Jh_n}
\def \OVn{\OVh_n}
\def \OVh{\OV^h}
\def \OV{{\mathcal O}(V)}

\def \refeq#1{equation (\ref{#1})}
\def \reff#1{(\ref{#1})}
\def \ra{\rightarrow}
\def \I{\mbox{${\mathcal I}$}}
\def \hom{\mbox{\rm Hom}}
\def \ie{\hbox{\it i.e.}}
\def \gl{\mbox{$\mathfrak{gl}$}}
\def \ord{\operatorname{ord}}
\def \tns{\otimes}
\def \mtns{\tns\cdots\tns}
\def \mplus{+\cdots+}
\def \mcom{,\cdots,}
\def \ms{\m_\star}
\def \h{h}
\def \Ah{\A^\h}
\def \kh{\k[[\h]]}
\def \k{\mbox{$\mathbb F$}}

\def \defi#1{{\em#1}}
\def \t{\tns}
\def \d{\delta}
\def \X#1#2{X_{#1#2}}
\def \l{\lambda}
\def \p{\partial}
\def \Pb{\{\cdot\,,\cdot\}}
\def \pp{[\pi_1,\pi_2]}
\def \tt#1#2#3{\p^{#1}\!\t\!\p^{#2}\!\t\!\p^{#3}}
\def \type#1#2#3{(#1,#2,#3)}
\def \({\left(}
\def \){\right)}

\def \C{\mbox{$\mathbb C$}}
\def \Z{\mbox{$\mathbb Z$}}
\def \N{\mbox{$\mathbb N$}}

\author{Michael Penkava}
  \address{University of Wisconsin,
           Department of Mathematics,
           Eau Claire, WI 54702-4004}
  \email{penkavmr@uwec.edu}
\author{Pol Vanhaecke}
  \address{University of California,
           1015 Department of Mathematics,
           Davis, CA 95616-8633}
  \email{vanhaeck@math.ucdavis.edu}
  \address{Universit\'e des Sciences et Technolgies de Lille,
           U.F.R. de Math\'ematiques,
           59655 Villeneuve D'Ascq, France}
  \email{Pol.Vanhaecke@Univ-Lille1.fr}
\thanks{The research of the first author was partially funded by grants {}from the University
        of Wisconsin, Eau Claire.}
\title{Deformation quantization of polynomial Poisson algebras}
\begin{document}
  \nocite{bir,wei,hue,omy1,omy2,omy3,dl,cp,dri,dou,hum}
\begin{abstract}
  This paper discusses the notion of a deformation quantization for an
  arbitrary polynomial Poisson algebra $\A$. We examine the Hochschild
  cohomology group $H^3(\A)$ and find that if a deformation of $\A$
  exists it can be given by bidifferential operators. We then compute
  an explicit third order deformation quantization of $\A$ and show
  that it comes {}from a quantized enveloping algebra. We show that
  the deformation extends to a fourth order deformation if and only if
  the quantized enveloping algebra gives a fourth order deformation;
  moreover we give an example where the deformation does not extend. A
  correction term to the third order quantization given by the
  enveloping algebra is computed, which precisely cancels the
  obstruction.
\end{abstract}

  \maketitle
  \tableofcontents
\section{Introduction}
Deformation theory for associative commutative algebras was first
considered by Gerstenhaber in \cite{ger}. A formal deformation of an
associative commutative algebra $\A$ over a ground field $\k$ is by
definition an associative multiplication $\star$ on $\Ah=\A[[h]]$,
\begin{equation}\label{defint}
  p\star q=pq+\h\m_1(p,q)+\h^2\m_2(p,q)+\cdots,
\end{equation}
where $pq$ denotes the original product of elements $p,\,q\in\A$. The
main tools which are used by Gerstenhaber are the Hochschild
cohomology groups $H^n(\A)$ (introduced in \cite{hoch}) and the
Gerstenhaber bracket $[\cdot\,,\cdot]$ (introduced in
\cite{gers}). In fact, the deformed product
\begin{equation*}
  \ms:\Ah\times\Ah\to\Ah
\end{equation*}
defines an associative product if and only if $[\ms,\ms]=0$, an
equation which can be rewritten by using the Hochschild coboundary
operator $\d$ as
\begin{equation*}
  \d\m_k=\frac12\sum_{i+j=k}[\m_i,\m_j]\qquad k=1,2,\dots
\end{equation*}
It is a fundamental fact that the right hand side of this equation is
a Hochschild 3-cocycle: if a deformation is associative up to order
$n$ then it extends to order $n+1$ if and only if some given 3-cocycle
is a coboundary, hence this question is cohomological in nature. One
immediate consequence is that the vanishing of $H^3(\A)$ implies that
every $n$-th order deformation extends to a formal deformation.

The relevance of deformation theory to physics was first pointed out
in \cite{bffls}. The main idea is that the non-commutative
(associative) operator product which appears in quantum mechanics is
a deformation of the commutative product of classical observables,
making deformation theory a tool for describing the transition {}from
classical to quantum mechanics. In this context one often speaks of a
deformation quantization (in the present paper we will reserve the
term deformation quantizations for deformations that alternate in the
sense that $\m_i$ is even or odd according to the parity of $i$). A
new object which appears in their treatment of deformation theory is a
Poisson bracket. Indeed, it is easy to show that if \reff{defint}
defines a deformation (of at least order 2) and $\pi_1$ is antisymmetric then
$\pi_1$ is a Poisson bracket.

This observation leads to the following question. Let $\A$ an
associative commutative algebra, equipped with a Poisson bracket, \ie,
an antisymmetric biderivation $\Pb:\A\times\A\to\A$ which satisfies
the Jacobi identity
\begin{equation}
  \{p,\{q,r\}\}+\{q,\{r,p\}\}+\{r,\{p,q\}\}=0, \qquad\hbox{for all }p,q,r\in\A.
\end{equation}
Does there exist a formal deformation \reff{defint} for which
$\m_1$ is given by the Poisson bracket,
$\m_1=\frac12\Pb$?
In two important cases an affirmative answer is given in \cite{bffls}:
when $\A$ is the Poisson algebra of functions on the
dual of a Lie algebra and when $\A$ is the algebra of functions on a
Poisson manifold which admits a flat connection.
In the  Lie algebra case, this result had already been observed by
Berezin (see \cite{ber}), who also pointed out the relation with the
enveloping algebra. The case of symplectic manifolds was settled later
by De Wilde and Lecomte (see \cite{dl}; for a geometric proof of
their result, see \cite{fed}).

The starting point of our research was to try to understand the case
of a general Poisson algebra by first investigating the case of
polynomial Poisson algebras (in any number of variables). For the
latter we use a subcomplex of the Hochschild complex, consisting of
differential operators. Indeed, for a polynomial ring, the third
cohomology group of the subcomplex of  differential operators maps
injectively
to the ordinary Hochschild cohomology group $H^3(\A)$.  Since the
Poisson bracket is a differential operator and the
Gerstenhaber bracket of differential operators is again a differential
operator, the problem of extending a given deformation quantization can
be described by the cohomology of differential operators,
in which one can do explicit computations more easily.
The fact that one can replace the Hochschild cohomology
groups with the smaller groups is not trivial and is based upon a
careful investigation of $H^3(\A)$. For a 3-cocycle $\ph$ we will show
that its flip symmetric part $\ph_+$, which is defined by
$\ph_+(p,q,r)=\frac12(\ph(p,q,r)+\ph(r,q,p))$ is always a coboundary
and that its flip antisymmetric part $\ph_-=\ph-\ph_+$ is a coboundary
if and only if
$\ph_-(\ph_-(p,q),r)+\ph_-(\ph_-(q,r),p)+\ph_-(\ph_-(r,p),q)=0$. In
either case we give a recursion formula for the cochain whose
coboundary $\ph$ is given and find that this cochain is a
bidifferential operator when $\ph$ is a tridifferential operator. A
characterization of $H^3(\A)$ can already be found in \cite{omy3}, but our
proofs have the advantage of being purely cohomological and allow the
latter conclusion.

Armed with the above explicit description of $H^3(\A)$ and explicit
formulas for the Gerstenhaber bracket and the Hochschild coboundary
operator we easily find that a deformation quantization of order three
always exists; moreover, using the Jacobi identity we can actually
write down an explicit formula for $\pi_2$. We also give an
explicit formula for $\pi_3$, as a result of a lot of non-trivial
computations which not only involve the Jacobi identity but also
its derivative. Thus we find that every polynomial Poisson algebra
admits a third order deformation quantization. Surprisingly enough
this seemingly ``natural'' deformation quantization does \emph{not}
(in general)
extend to a fourth order deformation.

This fact is even more striking once one realizes that the third order
deformation which we construct comes {}from a quantized universal
enveloping algebra, making this deformation most natural. We define
this enveloping algebra for any polynomial Poisson algebra $(\A,\Pb)$
as follows. First notice that $\A$ can be seen as the symmetric
algebra $\SV$ over a vector space $V$; then the Poisson bracket is a
linear map $\SV\bigotimes\SV\to\SV$. We take the tensor algebra $\TV$ of
$V$ and we consider  the two-sided ideal
$J^\h$ of  $\TVh$ generated by all elements of the form $x\t y-y\t
x-h\sg\{x,y\}$, where $x,y\in V$ and $\sg:\SV\to\TV$  is
the symmetrization map, defined by
\begin{equation*}
  \sg\left(\prod_{i=1}^n a_i\right)={\frac 1 {n!}}\sum_{p\in S_n}
      a_{p(1)}\t a_{p(2)}\t\cdots\t a_{p(n)}.
\end{equation*}
The quantized universal enveloping algebra is defined as
$\Uh=\TVh/J^h$. Notice that in the case of a linear Poisson
bracket we recover the usual definition of the enveloping algebra of a
Lie algebra. It is a well-known but non-trivial fact that for a linear
Poisson bracket
the enveloping algebra does give a deformation quantization in the
following way: the natural map $\SVh\to\Uh$ is a linear isomorphism,
so the product on $\Uh$ determines a product on $\SVh$,
which is a deformation quantization. In general, \ie, for non-linear
Poisson brackets the map $\SVh\to\Uh$ fails to be injective, but
surprisingly enough, for a general Poisson bracket it is injective
precisely up to order 3 (in $\h$). In fact, one computes an obstruction to the
injectivity of the map, which turns out to
coincide with the obstruction which we found earlier when trying to
extend the deformation to a fourth order deformation quantization.
Thus the third order deformation which we construct using Hochschild cohomology
extends to a fourth order deformation
precisely when the quantized enveloping algebra gives a fourth order
deformation. An explanation of this will be given in the text.

The recent result by Kontsevich, which states that every Poisson
manifold has a deformation quantization (see \cite{kon4}) was the
motivation for us to look what was ``wrong'' with our third
deformation. It is easy to see that one can always add any
antisymmetric biderivation to $\m_3$ and obtain a new, non-equivalent
(third order) deformation quantization. Even more, for some choices of
biderivation the third order deformation extends while for others it
doesn't. We will give such a biderivation for which the extension to a
fourth order deformation is always possible.  The check that it does
depends on a skillful use of the Jacobi identity, and the first and
second derivatives of the Jacobi identity.

The structure of this paper is as follows. In Section~2 we explain the
precise relation between deformation theory and Hochschild cohomology
for associative commutative algebras and we study the second
and third Hochschild cohomology groups for a polynomial algebra in
Sections~3 and~4. An explicit third order deformation for any
polynomial Poisson algebra is computed in Section~5 and in Section~6
we compute the obstruction for this deformation to extend to a fourth
order deformation. In Section~7 we introduce the quantized enveloping
algebra of a polynomial Poisson algebra and we show that our third
order deformation comes from this algebra. We show in Section~8 how to
modify the third order deformation such that it extends to a
fourth order deformation. In the final section a few examples with
very different characteristics are worked out, in particular we give
an example which shows that the third order deformation which is given
by the quantized universal enveloping algebra does not extend in
general.

\emph{Acknowledgements.}\quad The authors would like to thank
Alexander Astashkevich, Dmitry Fuchs, Josef Mattes, Bruno Nachtergaele
and Alan Weinstein for useful conversations. The first author would
also like to thank the mathematics department at the University of
California, Davis for providing office space during his two trips to
Davis to work on this project.
\section{Deformations of polynomial algebras}
In this section we briefly discuss the notion of deformation of a
commutative associative algebra $\A$ over a field $\k$. (We assume
throughout this paper that the characteristic of $\k$ is not 2.)
We describe the  obstruction to the existence of a deformation using the
Hochschild cohomology group $H^3(\A)$,
which has an explicit description in the case of a polynomial
algebra.

We will denote the product $pq$ of elements $p$, $q$ in $\A$ by
$\m(p,q)$.
Let $\h$ be a formal parameter and let $\Ah$ (resp.\ $\kh$) denote the
algebra of formal power series with coefficients in $\A$ (resp.\ in $\k$).
For $n\in\N$ we will also use the algebra $\An$ which is obtained {}from
$\Ah$ by dividing out by the ideal generated by $\h^{n+1}$. For elements $p,\,q\in
\Ah$ we write $p=q\mod \h^{n+1}$ when they project to the same elements in
$\An$.
\begin{dfn}\label{deformation}
  {An $\kh$-bilinear map
  \begin{equation*}
    \ms:\Ah\times\Ah\to\Ah
  \end{equation*}
  is called a \emph{ (formal) deformation} of $\A$ when it satisfies
  the associativity condition
  \begin{equation*}
    \ms(\ms(p,q),r)=\ms(p,\ms(q,r))
  \end{equation*}
  for all $p,\,q$ and $r$ in $\Ah$ and reduces to $\m$ on
  $\A\cong\Azero$, \ie,
  \begin{equation*}
    \ms(p,q)=\m(p,q)\mod \h.
  \end{equation*}

  More generally, when associativity merely holds on $\An$
  we say that $\ms$ defines an $n$-th order deformation. A first-order
  deformation is also called an \emph{infinitesimal deformation.}

  When a (formal) deformation has the additional property that for any
  $p,q\in\A$ the product $q\star p$ is obtained {}from $p\star q$ by
  applying the involution of $\Ah$ determined by $h\mapsto -h$, then
  we say that it defines a (formal) \emph{ deformation quantization}
  of $\A.$

  Two ($n$-th order or formal) deformations $\ms$ and $\ms'$ are
  called equivalent if there exists an $\kh$-linear map $F:\Ah\to\Ah$
  such that
  \begin{equation*}
    F(p)=p\mod\h
  \end{equation*}
  for any $p\in\A$ and such that $F(\ms(p,q))=\ms'(F(p),F(q))$ for any
  $p,\,q\in\A$.}
\end{dfn}

For $n\geq0$ the space of $n$-cochains is given by
\begin{equation*}
  C^n(\A)=\hom(\A^{n},\A),
\end{equation*}
and note that $\ms$ is given by a sequence of elements $\m_i$ in
$C^2(\A)$,
\begin{equation*}
  \ms=\m+\h\m_1+\h^2\m_2+\cdots.
\end{equation*}
A deformation $\ms$ is a
deformation quantization if $\pi_k$ is symmetric
when $k$ is even and antisymmetric when $k$ is odd.

The condition that $\ms$ be associative is most conveniently expressed
in cohomological language.
A graded bracket on cochains, called the \defi{Gerstenhaber bracket}
(see~\cite{ger}), is given by
\begin{multline*}
  [\ph,\psi](p_1\mcom p_{m+n-1})=\\
    \sum_{k=1}^m(-1)^{(k-1)(n-1)}\ph(p_1\mcom p_{k-1},
  \psi(p_k \mcom p_{k+n-1}), p_{k+n}\mcom p_{m+n-1})\\
    -(-1)^{(m-1)(n-1)}\times\\
  \sum_{k=1}^n(-1)^{(k-1)(m-1)}\psi(p_1\mcom p_{k-1},
  \ph(p_k \mcom p_{k+m-1}), p_{k+m}\mcom p_{m+n-1}),
\end{multline*}
for $\ph\in C^m(\A),\,\psi\in C^n(\A)$. We will use the bracket in the
case of elements $\ph,\,\psi\in C^2(\A)$ in which case
$[\ph,\psi]=[\psi,\ph]$ and the formula specializes to
\begin{equation*}
  [\ph,\psi](p,q,r)=\ph(\psi(p,q),r)-\ph(p,\psi(q,r))
  +\psi(\ph(p,q),r)-\psi(p,\ph(q,r)).
\end{equation*}
In view of the above formula the proof of the following lemma is
trivial.
\begin{lma}\label{ass}
  An element $\ph\in C^2(\A)$ defines an associative multiplication on
  $\A$ if and only if $[\ph,\ph]=0$.
\end{lma}
It follows that $\ms$ is associative if and only if
\begin{equation}\label{cohom_0}
  0=[\ms,\ms]=2\h[\m,\pi_1]+\h^2(2[\m,\pi_2]+[\pi_1,\pi_1])+\cdots
\end{equation}
The cochains form a complex $C^\bullet(\A)$ for the Hochschild coboundary
operator
\begin{equation*}
 \d:C^n(\A)\to C^{n+1}(\A)
\end{equation*}
which is defined  by
\begin{multline*}
  \d\ph(p_1\mcom p_{n+1})=p_1\ph(p_2\mcom p_{n+1})\\
  +\sum_{k=1}^n(-1)^k\ph(p_1\mcom p_{k-1},p_kp_{k+1},p_{k+2}\mcom p_{n+1})
  +(-1)^{n+1}\ph(p_1\mcom p_n)p_{n+1},
\end{multline*}
for $\ph\in C^n(\A)$ (see \cite{hoch}).  The $n$-th cohomology group
of this complex will be denoted by $H^n(\A)$.  It is easy to see that
the Hochschild coboundary operator $\delta$ can be written in terms of
the Gerstenhaber bracket as
\begin{equation*}
  \d\ph=-[\ph,\m],
\end{equation*}
so that the associativity condition (\ref{cohom_0}) can be expressed
by an infinite list of relations
\begin{align}\label{cohom}
  &\d\m_1=0,\notag\\
  &\d\m_2=\tfrac12[\m_1,\m_1],\notag\\
  &\d\m_3=[\m_1,\m_2],\\
  &\d\m_4=[\m_1,\m_3]+\tfrac12[\m_2,\m_2],\notag\\
  &\phantom{dddddd}\vdots\notag
\end{align}
More precisely, if the cochains $\m_1,\dots,\m_{n-1}$ define an
$(n-1)$-th order deformation of $\m$ then it extends to an $n$-th order
deformation if and only if the equation
\begin{equation}\label{massey}
  \d\m_n={\tfrac12}\sum_{i+j=n}[\m_i,\m_j]
\end{equation}
has a solution $\m_n$. If such a solution exists it is clearly unique up to
addition of any cocycle $\d\ph,\,\ph\in C^1(\A)$. As for its existence
it is important to note that the right hand side in (\ref{massey}) is always
a cocycle:
\begin{align*}
  \sum_{i+j=n}\d[\m_i,\m_j]
  &=-\sum_{i+j=n}[[\m_i,\m_j],\m],\\
  &= 2\sum_{i+j=n}[[\m_j,\m],\m_i],\\
  &=-\sum_{i+k+l=n}[[\m_k,\m_l],\m_i],\\
  &=0.
\end{align*}
In this computation the Jacobi identity
\begin{equation}\label{big_Jacobi}
  [[\ph,\chi],\psi]+[[\chi,\psi],\ph]+[[\psi,\ph],\chi]=0
\end{equation}
which is valid for any 2-cochains $\ph,\,\psi$ and $\chi$, was used
twice (when the characteristic of $\k$ is 3 then the equation
$[[\ph,\ph],\ph]=0$, which holds for any 2-cochain $\ph$, but
\emph{not} as a consequence of (\ref{big_Jacobi}), is also used.)
The upshot is that the extendibility of a deformation of order $n-1$
to a deformation of order $n$ depends on whether or not a certain
Hochschild 3-cocycle is a coboundary. However, the particular
$\m_{n}$ chosen for the extension of the deformation to order $n$ will
have a pronounced impact on the further extendibility of the
deformation. If the deformation does not extend to order $n+1$, it may
be that a different choice of $\m_{n}$ would allow such an
extension. Moreover, this effect is not limited to the next term in
the extension, so that the extendibility of an extension up to order
$n$ is influenced by all of the choices of the cochains $\m_k$ for
$k<n$.

One can describe the extendibility of the deformation in terms of
Massey powers of $\m_1$(see \cite{mass:def}), so that there is an
extension of order $n$ when the $n$-th Massey power of $\m_1$
vanishes, but this description does not yield any immediate
advantage, since the problem of computation of the Massey powers
may be more difficult to solve than the problem of finding an
explicit sequence of cochains yielding a deformation.

Another important consideration is the uniqueness of deformations, which
is partially governed by the second Hochschild cohomology group
$H^2(\A)$. The following lemma shows that if two deformations differ by
a coboundary then they are equivalent.
\begin{prp}\label{isom}
  If the $n$-th cochain $\m_n$ in a formal (resp.\ $m$-th order with
  $n\le m$) deformation $\sum\h^i\m_i$ is altered by a coboundary then
  the new $n$-th order deformation extends to an equivalent formal
  (resp.\ $m$-th order) deformation.
\end{prp}
\begin{proof}
  In the case of a formal deformation, let us denote the coboundary
  which is added to $\pi_n$ by $\d E$, where $E\in C^1(\A)$.  Define an
  $\kh$-linear map $F:\Ah\to\Ah$ by $$F(p)= p+\h^nE(p) $$ with inverse
  $$F^{-1}(p)=p-\sum\h^{kn}E^k(p) $$ for any $p\in\A$. Then
  $\ms'(p,q):=F^{-1}(\ms(F(p),F(q)))$ is a new (equivalent)
  deformation whose first $n$ cocycles $\m_i'$ coincide with the
  cocycles $\m_i$ and $\m_n'=\m_n+\d E$.
\end{proof}

\section{Hochschild cohomology}\label{hoch}
Our aim in this section is to analyze $H^2(\A)$ and $H^3(\A)$ more
thoroughly. We will give an explicit characterization in the case of
polynomial algebras, in Theorems \ref{2sym}, \ref{3skew} and
\ref{3sym}.  These results should be regarded as classical, and
complete proofs of these theorems appear in \cite{omy3}; but our
treatment here is more straightforward, relying on simple
cohomological arguments.  In the proofs below, for a
polynomial algebra given by an ordered basis (free generating set)
$\{x_i\}_{i\in\I}$, we will denote elements of the basis by the
letters $x$ and $y$, while arbitrary polynomials will be denoted by
the letters $p$, $q$, $r$ and $s$. For $x\in\I$ the statement $x\le p$
means that the basis elements appearing in the monomials in $p$ have
index greater than or equal to that of $x$, so that in particular
$x\le c$ for any constant $c$.

Any 2-cochain $\ph$ can be uniquely decomposed as the sum of a
symmetric cochain $\ph^+$ and an antisymmetric cochain $\ph^-$. Then $\ph$
is a cocycle precisely when both its symmetric and antisymmetric parts
are cocycles.  To see this fact, suppose that $\ph$ is a
2-cocycle. Then
$\delta\ph(p,q,r)=p\ph(q,r)-\ph(pq,r)+\ph(p,qr)-\ph(p,q)r$. Let
$\bar\ph(p,q)=\ph(q,p)$. Then
\begin{equation*}
  \delta\bar\ph(p,q,r)=p\ph(r,q)-\ph(r,pq)+\ph(qr,p)-
  \ph(q,p)r=-\delta\ph(r,q,p)=0.
\end{equation*}
Since $\ph^+$ and $\ph^-$ are linear combinations of $\ph$ and
$\bar\ph$, this shows the desired result.  Furthermore, the coboundary
of any 1-cochain is symmetric, which is immediate {}from the fact that
if $\lambda$ is a 1-cochain, then
\begin{equation*}
  \delta\lambda(p,q)=p\lambda(q)-\lambda(pq)+\lambda(p)q.
\end{equation*}
This implies that each antisymmetric 2-cocycle determines a distinct
cohomology class.  Any biderivation is a 2-cocycle, since for a
biderivation $\ph$
\begin{equation*}
  \delta\ph(a,b,c)=a\ph(b,c)-\ph(ab,c)+\ph(a,bc)-\ph(a,b)c=-\ph(a,c)b+b\ph(a,c)=0.
\end{equation*}
Furthermore, any antisymmetric 2-cochain is a cocycle precisely when
it is a bideri\-vation.  To see this, note that if $\ph$ is an
antisymmetric cocycle, then
\begin{equation*}
  \delta\ph(a,b,c)-\delta\ph(c,a,b)+\delta\ph(b,c,a)=2(\ph(a,bc)-b\ph(a,c)-
  \ph(a,b)c).
\end{equation*}
Since the left hand side vanishes, $\ph$ is a biderivation.  These
remarks hold for an arbitrary commutative algebra $\A$, but when $\A$
is a polynomial algebra, we have a more complete characterization of
$H^2(\A)$.
\begin{thm}\label{2sym}
  Let $\A$ be a polynomial algebra. Then a 2-cocycle $\ph$ is a
  coboundary precisely when it is symmetric. Furthermore, the cochain
  $\lambda$ satisfying $\delta\lambda=\ph$ can be chosen arbitrarily
  for basis elements. In particular, it can be chosen to satisfy
  $\lambda(x)=0$ when $x$ is a basis element.
\end{thm}
\begin{proof}
  For a symmetric 2-cocycle $\ph$, we construct recursively a
  1-cochain $\lambda$ whose coboundary coincides with $\ph$.  Now
  \begin{equation*}
    \delta\ph(1,1,q)=\ph(1,q)-\ph(1,q)+\ph(1,q)-\ph(1,1)q=0,
  \end{equation*}
  so that $\ph(1,q)=\ph(1,1)q$. Let $\lambda(1)=\ph(1,1)$, and define
  $\lambda(x)$ arbitrarily when $x$ has degree 1.  The property
  $\d\lambda=\ph$ holds precisely when
  \begin{equation*}
    \lambda(pq)=p\lambda(q)+q\lambda(p)-\ph(p,q).
  \end{equation*}
  When either $p$ or $q$ is constant, this equation holds by the
  preceding remarks; otherwise $\lambda$ is evaluated at terms of
  lower degree on the right hand side, so the left hand side is
  defined recursively by this formula. But we need to check that if
  $pq=p'q'$ then the right hand sides of the decomposition above
  agree.
%
Equivalently, it is enough to check that expanding $\lambda(pqr)$ does not
depend on whether $p$ or $q$ is factored out first.
  Consider the expansion
  \begin{align*}
    \lambda(pqr)&=p\lambda(qr)+qr\lambda(p)-\ph(p,qr)\\
    &=p(q\lambda(r)+r\lambda(q)-\ph(q,r))+qr\lambda(p)-\ph(p,qr).
  \end{align*}
  Similarly,
  \begin{equation*}
    \lambda(qpr)=q(p\lambda(r)+r\lambda(p)-\ph(p,r))+pr\lambda(q)-\ph(q,pr).
  \end{equation*}
  Using the symmetry of $\ph$, we find that the difference of these
  two expressions is $\d\ph(p,r,q)$, and so vanishes.
\end{proof}

Let us now turn our attention to the third Hochschild cohomology
group, wherein lie the obstructions to extension of a deformation to
higher order.  A 3-cochain $\psi$ is called flip symmetric if it
satisfies $\psi(p,q,r)=\psi(r,q,p)$, and flip antisymmetric if
$\psi(p,q,r)=-\psi(r,q,p)$. Every 3-cochain $\psi$ can be uniquely
decomposed as the sum of a flip symmetric cochain $\psi^+$ and a flip
antisymmetric cochain $\psi^-$.  Moreover, $\psi$ is a cocycle
precisely when $\psi^+$ and $\psi^-$ are cocycles. Furthermore, the
coboundary of a symmetric 2-cochain is flip antisymmetric, and the
coboundary of an antisymmetric 2-cochain is flip symmetric. The
\emph{Jacobi map} $J:C^3(\A)\to C^3(\A)$ is given by
\begin{equation*}
  J\psi(p,q,r)=\psi(p,q,r)+\psi(q,r,p)+\psi(r,p,q).
\end{equation*}
Then, if $\psi$ is the coboundary of a symmetric cochain, it satisfies
the {\em Jacobi identity} $J\psi=0$. For any 3-cocycle $\psi$,
$\psi(1,1,1)=\delta\psi(1,1,1,1)=0$ and
$\psi(1,p,1)=\delta\psi(1,1,p,1)=0$.  Suppose that
$\psi(p,1,1)=\psi(1,1,p)=0$ for all $p$.  Then
$\psi(1,p,q)=\delta\psi(1,1,p,q)$, $\psi(p,1,q)=\delta\psi(p,1,1,q)$,
and $\psi(p,q,1)=\delta\psi(p,q,1,1)$, so these terms vanish for all
$p$ and $q$.  These remarks are easy to check, and apply to any
commutative algebra $\A$, not just a polynomial algebra.
\begin{thm}\label{3skew}
  Let $\A$ be a polynomial algebra and suppose that $\psi$ is a
  3-cocycle. Then $\psi$ is flip symmetric if and only if $\psi$ is a
  Hochschild coboundary of an antisymmetric cochain $\ph$.
\end{thm}
\begin{proof}
  By the above remarks we only need to verify that a flip symmetric
  cocycle $\psi$ is a coboundary of an antisymmetric cochain. If we
  define $\theta(p,1)=\psi(1,1,p)=-\theta(1,p)$, and extend $\theta$
  in an arbitrary manner to an antisymmetric cochain, then
  $\delta\theta(1,1,p)=-\psi(1,1,p)$, so that by replacing $\psi$ by
  $\psi+\delta\theta$, we may assume that $\psi(1,1,p)=0$, so that
  $\psi$ vanishes when any of its arguments is a constant.  We define
  $\ph$ recursively, by first setting
  $\ph(1,1)=\ph(1,p)=\ph(p,1)=0$. In addition, let us assume that
  $\ph(x,y)$ is defined in an arbitrary manner for basis elements $x$
  and $y$.  Consider the following equalities which must be satisfied
  if $\psi=\delta\ph$.
  \begin{align*}
    \psi(p,q,r)&=p\ph(q,r)-\ph(pq,r)+\ph(p,qr)-\ph(p,q)r,\\
    \psi(p,r,q)&=p\ph(r,q)-\ph(pr,q)+\ph(p,rq)-\ph(p,r)q,\\
    \psi(r,p,q)&=r\ph(p,q)-\ph(rp,q)+\ph(r,pq)-\ph(r,p)q.
  \end{align*}
  Adding the first and third and subtracting the second of these
  equations, and using the desired antisymmetry property for $\ph$
  yields the following equation:
  \begin{equation*}
    2\ph(pq,r)=2p\ph(q,r)+2q\ph(p,r)-\psi(p,q,r)+\psi(p,r,q)-\psi(r,p,q),
  \end{equation*}
  The expression above is evidently symmetric in $p$ and $q$, and
  holds when either $p$ or $q$ is constant, because $\psi$ vanishes
  when any of its arguments is a constant.  Otherwise the right hand
  side involves terms of smaller degree than the left, so we obtain a
  recursive definition of $\ph$.  To check that this definition does
  not depend on the decomposition of the term $pq$, one should check
  that expanding a term of the form $\ph(pqr,s)$ in two ways leads to
  the same result. Factoring $pqr$ as $p$ times $qr$ one obtains
  \begin{align*}
    2\ph(pqr,s)&=2p\ph(qr,s)+2qr\ph(p,s)-\psi(p,qr,s)+\psi(p,s,qr)-\psi(s,p,qr)\\
    &=p(2q\ph(r,s) +2r\ph(q,s)-\psi(q,r,s)+\psi(q,s,r)-\psi(s,q,r))+\\
    &\qquad 2qr\ph(p,s) -\psi(p,qr,s)+\psi(p,s,qr)-\psi(s,p,qr).
  \end{align*}
  Subtracting the expression obtained by interchanging the roles of
  $p$ and $q$, the terms involving $\ph$ drop out and we are left with
  \begin{align*}
    &-p\psi(q,r,s)+p\psi(q,s,r)-p\psi(s,q,r)
    -\psi(p,qr,s)+\psi(p,s,qr)-\psi(s,p,qr)\\
    &+q\psi(p,r,s)-q\psi(p,s,r)+q\psi(s,p,r)
    +\psi(q,pr,s)-\psi(q,s,pr)+\psi(s,q,pr)\\ &=
    \delta\psi(s,p,r,q)-\delta\psi(p,r,q,s)-
    \delta\psi(p,s,r,q)+\delta\psi(q,s,r,p),
  \end{align*}
  which is zero.  To verify that $\ph$ is antisymmetric, we need to
  compute $\ph(pq,rs)+\ph(rs,pq)$, using antisymmetry of lower degree
  terms.  We have
  \begin{align*}
    2\ph(pq,rs)&=2p\ph(q,rs)+2q\ph(p,rs)-\psi(p,q,rs)+\psi(p,rs,q)-\psi(rs,p,q)\\
    &=p(2r\ph(q,s)+2s\ph(q,r)+\psi(r,s,q)-\psi(r,q,s)+\psi(q,r,s))+\\
    &\qquad q(2r\ph(p,s)+2s\ph(p,r)+\psi(r,s,p)-\psi(r,p,s)+\psi(p,r,s))+\\
    &\qquad \psi(p,q,rs)+\psi(p,rs,q)-\psi(rs,p,q).
  \end{align*}
  {}From this we see that $2(\ph(pq,rs)+\ph(rs,pq))$ equals
  \begin{multline*}
    \delta\psi(p,q,r,s)-\delta\psi(p,r,q,s)+\delta\psi(q,s,r,p)\\
    +\delta\psi(s,q,p,r)+ \delta\psi(r,s,p,q)-\delta\psi(q,s,p,r),
  \end{multline*}
  and thus vanishes. It is immediately checked that $\delta\ph=\psi$.
\end{proof}
\begin{thm}\label{3sym}
  Let $\A$ be a polynomial algebra and suppose that $\psi$ is a
  3-cocycle.  Then $\psi$ is a Hochschild coboundary of a symmetric
  cochain $\ph$ if and only if $\psi$ is flip antisymmetric and
  satisfies the Jacobi identity $J\psi=0.$ Moreover, if we take an
  ordered basis of $\A$, then $\ph$ can be chosen to satisfy
  $\ph(x,p)=0$ whenever $x$ is a basis element satisfying $x\le p$.
\end{thm}
\begin{proof}
  As in the previous theorem, we reduce to the case where
  $\psi(p,1,1)=0$. Take an ordered basis of $\A$.  Define $\ph$ by
  $\ph(1,p)=0$ for all $p$ and $\ph(x,p)=0$ when $x\le p$.  We extend
  the definition recursively by setting
  \begin{equation*}
    \ph(xp,q)=x\ph(p,q)+\psi(q,p,x)=\ph(q,xp),
  \end{equation*}
  when $x\le q$ and $x\le p$. To show $\ph$ is well defined and
  symmetric, we only need to show that if $x$ is a basis element
  satisfying $x\le p$ and $x\le q$, then the expansion of $\ph(xp,xq)$
  yields the same result as the expansion of $\ph(xq,xp)$.  Now
  \begin{equation*}
    \ph(xp,xq)=x\ph(p,xq)+\psi(xq,p,x)=x(x\ph(q,p)+\psi(p,q,x))+\psi(xq,p,x),
  \end{equation*}
  so that $\ph(xp,xq)-\ph(xq,xp)=\delta\psi(x,p,q,x)=0$.

  To show that $\delta\ph=\psi$, we note that if any of $p$, $q$ or
  $r$ is constant, then both $\psi(p,q,r)$ and $\delta\ph(p,q,r)$
  vanish (for the vanishing of the former, see the remarks preceeding
  Theorem \ref{3skew}). We may proceed by induction on the sum of the
  degrees of $p$, $q$ and $r$.  If $p$ can be factored as $xp'$, where
  $x$ satisfies $x\le p'$, $x\le q$ and $x\le r$, then
  \begin{align*}
    \delta\ph(xp',q,r)&=xp'\ph(q,r)-\ph(xp'q,r)+\ph(xp',qr)-\ph(xp',q)r\\
    &=xp'\ph(q,r) -x\ph(p'q,r)-\psi(r,p'q,x)\\
    &\qquad+x\ph(p',qr)+\psi(qr,p',x) -rx\ph(p',q)-r\psi(q,p',x)\\
    &=x(p'\ph(q,r)-\ph(p'q,r)+\ph(p',qr)-\ph(p',q)r)\\
    &\qquad-\psi(r,p'q,x)+\psi(qr,p',x)-r\psi(q,p',x)\\
    &=x\psi(p',q,r)
    -\psi(r,p'q,x)+\psi(qr,p',x)-r\psi(q,p',x)\\&=\psi(xp',q,r).
  \end{align*}
  On the other hand, if we can express $r=xr'$, where $x\le p$, $x\le
  q$ and $x\le r$, then
  \begin{equation*}
    \psi(p,q,xr')=-\psi(xr',q,p)=
    -\delta\ph(xr',q,p)=\delta\ph(p,q,xr'),
  \end{equation*}
  since $\psi$ is flip antisymmetric, and the coboundary of any
  symmetric cochain is also flip antisymmetric.  The only other
  possibility is that $q=xq'$, where $x\le q'$, $x\le p$ and $x\le
  r$. But then we have
  \begin{align*}
    \psi(p,xq',r)&=-\psi(xq',r,p)-\psi(r,p,xq')\\
    &=-\delta\ph(xq',r,p)-\delta\ph(r,p,xq')\\ &=\delta\ph(p,xq',r),
  \end{align*}
  using the Jacobi identity $J\psi=0$ and the fact that the coboundary
  of any symmetric cochain satisfies the Jacobi identity. Note that it
  is only at this last step that the Jacobi identity is used.
\end{proof}
For simplicity in the proof above, we constructed $\ph$ so that
$\ph(x,p)=0$ for a basis element $x$ satisfying $x\le p$. But for a
polynomial algebra, one can always define $\ph(x,p)$ for $x\le p$ in
an arbitrary manner and extend the definition to a cocycle, as we show
below.  Thus we could have assumed in Theorem \ref{3sym} that
$\ph(x,p)$ is defined arbitrarily for $x\le p$.
\begin{prp}
  Let $\A$ be a polynomial algebra with an ordered basis, and suppose
  that $\ph(x,p)$ is any cochain defined for $x\le p$, satisfying
  $\ph(x,1)=0.$ Then $\ph$ extends uniquely to a symmetric cocycle
  satisfying $\ph(1,1)=0$.
\end{prp}
\begin{proof}
  {}From the condition $\delta\ph(x,p,q)=0$ one derives the property
  \begin{equation*}
    \ph(xp,q)=x\ph(p,q)+\ph(x,pq)-\ph(x,p)q.
  \end{equation*}
  If either $p$ or $q$ is constant, then the formula above holds
  trivially.  Otherwise, if $x\le p$ and $x\le q$, then the left hand
  side is defined recursively by the right hand side. The consistency
  and the symmetry condition $\ph(xp,xq)=\ph(xq,xp)$ follow {}from
  \begin{align*}
    \ph(xp,xq)&=x\ph(p,xq)+\ph(x,xpq)-\ph(x,p)xq\\
              &=x^2\ph(q,p)+x\ph(x,pq)-\ph(x,q)xp+\ph(x,xpq)-\ph(x,p)xq.
\end{align*}
  If $\ph(q,p)=\ph(p,q)$, then the above formula is already symmetric
  in $p$ and $q$, so the check of consistency and symmetry is trivial.
  To see that $\delta\ph(p,q,r)=0$, consider the case when $p=xp'$,
  where $x\le p'$, $x\le q$ and $x\le r$. Then
  \begin{align*}
    \delta\ph(xp',q,r)
    &=xp'\ph(q,r)-\ph(xp'q,r)+\ph(xp',qr)-\ph(xp',q)r\\
    &=xp'\ph(q,r)-x\ph(p'q,r)-\ph(x,p'qr)+\ph(x,p'q)r+x\ph(p',qr)+\\
    &\qquad\ph(x,p'qr)-\ph(x,p')qr-x\ph(p',q)r-\ph(x,p'q)r+\ph(x,p')qr\\
    &=x(p'\ph(q,r)-\ph(p'q,r)+\ph(p',qr)-\ph(p',q)r)\\
    &=x\d\ph(p,q,r),
  \end{align*}
  which is zero by the induction hypothesis.  The other cases follow
  {}from the flip antisymmetry and the Jacobi identity $J(\delta\ph)=0$.
\end{proof}
When applied to deformation theory Theorems \ref{3skew} and \ref{3sym}
lead to the following result.
\newpage
\begin{thm}\label{upto3}\quad\par
  \begin{enumerate}
    \item If $\m+\h\m_1$ is an infinitesimal deformation of $\A$ which
    extends to a second order deformation, then $\m_1^-$ is an
    antisymmetric biderivation which satisfies the usual Jacobi
    identity:
    \begin{equation}\label{Jacobi2}
      \m_1^-(\m_1^-(p,q),r)+\m_1^-(\m_1^-(q,r),p)+\m_1^-(\m_1^-(r,p),q)=0,
    \end{equation}
    so that $\m_1^-$ determines a Poisson algebra structure on $\A$.
    \item If $\A$ is a polynomial Poisson algebra, then the converse
    is true.  More precisely, if $\m+\h\m_1$ is an infinitesimal
    deformation such that $\m_1^-$ satisfies the usual Jacobi identity
    \reff{Jacobi2}, then the infinitesimal deformation extends to a
    second order deformation $\m+\h\m_1+\h^2\m_2$.  Furthermore, if
    $\pi_1$ is antisymmetric then $\m_2$ can be chosen to be symmetric
    in which case the deformation can be extended to order 3. In
    particular, any Poisson algebra structure on a polynomial algebra
    determines a deformation quantization of order 3.
  \end{enumerate}
\end{thm}

\begin{proof}
  If $\ph$ is antisymmetric, then $[\ph,\ph]$ satisfies the Jacobi
  identity $J([\ph,\ph])=0$ precisely when $\ph$ satisfies the usual
  Jacobi identity (\refeq{Jacobi2}).  With this remark, the statements
  in the theorem follow {}from our previous results.
\end{proof}
For a given Poisson algebra $(\A,\Pb)$ we will say that a deformation
$\ms=\m+\h \m_1+\h^2\m_2+\cdots$ of $\A$, in the sense of Definition
\ref{deformation}, defines a
\emph{deformation} of $(\A,\Pb)$ when $\m_1=\frac12\Pb$.

\section{Hochschild cohomology and differential operators}\label{hochdiff}
In this section we will assume that $\A$ is a polynomial algebra with
a fixed basis $\{x_i\}_{i\in\I}$ over a field $\k$ of characteristic
0.  We will give a characterization of Hochschild cochains in terms of
(possibly infinite order) differential operators.  First, let us
establish some conventions on our terminology.  For a basis element
$x_i$ of $\A$ we will denote the derivation $\p/\p x_i$ by $\p^i$. For
a multi-index $I=(i_1\mcom i_m)$, $\p^I$ will stand for the
differential operator $\p^{i_1}\dots\p^{i_m}$, $x_I$ will stand for
$x_{i_1}\dots x_{i_m}$, and $|I|=m$ is  its order.  For a polynomial $p$, we will
denote $\p^I(p)$ by $p^I$ and $\p^I(x_{I})$ by $I!$.  The multi-index $I$
is said to be \emph{non-decreasing} if $i_1\le\dots\le i_m$.  Also, we
shall write $I<I'$ to indicate that $I$ is obtained by removing some
of the indices in $I'$.  By $\p^{I_1}\mtns \p^{I_n}$ we shall denote
the $n$-differential operator of \emph{order} $|I_1|\mplus |I_n|$
given by
\begin{equation*}
  \p^{I_1}\mtns \p^{I_n}(p_1\mcom p_n)=\p^{I_1}(p_1)\dots \p^{I_n}(p_n).
\end{equation*}
The differential operator is said to have \emph{type} $(|I_1|\mcom
|I_n|)$. An expression of the form
\begin{equation*}
  \ph=\sum_{I_1\mcom I_n}\ph_{I_1\mcom I_n}\p^{I_1}\mtns \p^{I_n},
\end{equation*}
where $\ph_{I_1\mcom I_n}$ are polynomials, and we sum over all
non-decreasing multi-indices, gives a well-defined $n$-cochain on the
polynomial algebra. When only finitely many non-zero terms appear then
we say that $\ph$ is a \emph{(finite order) differential operator,}
otherwise such an expression is called a \emph{formal differential
operator.}  The \emph{order} of a differential operator $\ph$ is the
largest $m$ for which there is a nonzero term in $\ph$ of order $m$.
Every $n$-cochain can be expressed as a formal differential operator,
since we can solve for the polynomials $\ph_{I_1\mcom I_n}$ above
recursively by
\begin{equation*}
  I_1!\dots I_n! \ph_{I_1\mcom I_n}=\ph(x_{I_1}\mcom x_{I_n})-
  \hskip -.25in
  \sum_{(J_1\mcom J_n)<(I_1\mcom I_n)}
  \hskip -.25in
  \ph_{J_1\mcom J_n}\p^{J_1}(x_{I_1})\dots \p^{J_n}(x_{I_n}).
\end{equation*}
In the following lemma, which characterizes when an $n$-cochain is a
differential operator, we use the notation $(x-k)^I$ to stand for the
product $(x_1-k_1)^{i_1}\dots(x_s-k_s)^{i_s}$ for $k\in \k^{\I}$, and
$I=(i_1,\dots,i_s)$.
\begin{lma}\label{finiteorder}
  An $n$-cochain is a (finite order) differential operator precisely
  when there is some $N$ such that for any $k\in\k^{\I}$,
  \begin{equation*}
    \ph((x-k)^{I_1}\mcom (x-k)^{I_n})(k)=0,
  \end{equation*}
  whenever $|I_1|\mplus |I_n|\ge N$.
\end{lma}
\begin{proof}
  If $\ph$ is a (finite order) differential operator it suffices to
  take $N=1+\ord\ph$. On the other hand, if the order of $\ph$ is
  infinite, we may find for any $N$ a non-zero $\ph_{I_1\mcom I_n}$,
  with $|I_1|\mplus |I_n|\ge N$, in particular this polynomial is
  non-zero at some point $(k_1\mcom,k_n)$. Then
  \begin{equation*}
    \ph((x-k)^{I_1}\mcom (x-k)^{I_n})(k)=I_1!\dots I_n!\ph_{I_1\mcom I_n}(k)\neq0.
  \end{equation*}
\end{proof}
In the proof of the above lemma, it was necessary to evaluate a
polynomial at a point. This is the only place where the arguments in
this section cannot be extended to the ring of formal power series in
the variables $\{x_i\}_i\in\I$, because evaluation at a point is not
well defined. By examining the recursion formulas in Theorems
\ref{3skew} and \ref{3sym} and applying Lemma \ref{finiteorder}, one
sees that if $\psi$ is a differential operator, then the cochain $\ph$
satisfying $\delta\ph=\psi$ constructed in these theorems is also a
differential operator. The cochain $\lambda$ constructed in Theorem
\ref{2sym} will also be a differential operator when $\ph$ is a
bidifferential operator.

The notation $I'+I''=I$ will be used to indicate a partitioning of the
indices of the nondecreasing multi-index $I$ into two nondecreasing
multi-indices $I'$ and $I''$.  Then we obtain a very simple
description of the action of the Hochschild coboundary operator on
cochains which are given by differential operators, namely if $p\in\A$
and $J$ is any multi-index then
\begin{equation*}
  \d(p\p^I)=-\sum_{I'+I''=I}
  p\p^{I'}\tns \p^{I''}.
\end{equation*}
In general, if $\alpha$ is an $m$-cochain, and $\beta$ is an $n$-cochain,
then the $(m+n)$-cochain $\alpha\tns\beta$ is given by
\begin{equation*}
  \alpha\tns\beta(p_1\mcom p_{m+n})=
  \alpha(p_1\mcom p_m)\beta(p_{m+1}\mcom p_{m+n}).
\end{equation*}
The Hochschild coboundary operator acts as a graded derivation with
respect to this product, \ie,
\begin{equation*}
  \d(\alpha\tns\beta)=\d(\alpha)\tns \beta+\s{m}\alpha\tns\delta(\beta).
\end{equation*}
{}From this we obtain the following useful expression for the
coboundary of a bidifferential operator,
\begin{equation}\label{delta}
       \d(p\p^J\t\p^K)=p(\d\p^J)\t\p^K-p\p^J\t(\d\p^K).
\end{equation}
{}From the above, we see that the coboundary of an $n$-differential
operator of order $m$ is an $(n+1)$-differential operator of order
$m$.  The following theorem is an easy consequence of the above
remarks.
\begin{thm}\label{diffcob}
  Suppose that $\A$ is a polynomial algebra, and $\psi$ is an
  $n$-differential operator of order $m$.  If $\psi$ is a Hochschild
  coboundary, then we can find an $(n-1)$-differential operator $\ph$
  of order $m$ such that $\d\ph=\psi$.
\end{thm}
Let us denote
\begin{equation*}
  (p\p^{I_1}\mtns\p^{I_n})^J=\sum_{J_0\mplus J_n=J}
  p^{J_0}\p^{I_1J_1}\mtns \p^{I_nJ_n}.
\end{equation*}
Then the bracket of differential operators is given by
\begin{multline*}
  [p\p^{I_1}\mtns\p^{I_m},q\p^{J_1}\mtns \p^{J_n}]=\\
  \sum_{k=1}^m
  \s{(k-1)(n-1)}
  p\p^{I_1}\mtns\p^{I_{k-1}}\tns(q\p^{J_1}\mtns\p^{J_n})^{I_k}\tns
  \p^{I_{k+1}}\mtns\p^{I_m}\\
  -\s{(m-1)(n-1)}\times\\
  \sum_{k=1}^n
  \s{(k-1)(m-1)}
  q\p^{J_1}\mtns\p^{J_{k-1}}\tns(p\p^{I_1}\mtns\p^{I_m})^{J_k}\tns
  \p^{J_{k+1}}\mtns\p^{J_n}.
\end{multline*}
In particular, we obtain the following formula for the bracket of
2-cochains.
\begin{multline}\label{GHbracket}
  [p\p^{I_1}\t\p^{I_2},q\p^{J_1}\t\p^{J_2}]=\\
  p(q\p^{J_1}\tns\p^{J_2})^{I_1}\tns\p^{I_2}-
  p\p^{I_1}\tns(q\p^{J_1}\tns\p^{J_2})^{I_2}\\
  +q(p\p^{I_1}\tns\p^{I_2})^{J_1}\tns\p^{J_2}-
  q\p^{J_1}\tns(p\p^{J_1}\tns\p^{I_2})^{J_2}.
\end{multline}
The formula above will be used in the calculations of the
second and third order deformations of a polynomial Poisson algebra,
which will be carried out in the next section.

Finally, we apply the results of this section to show that any
deformation of a polynomial algebra is equivalent to one which is
given by differential operators.
\begin{thm}\label{diffquant}
  Any deformation (deformation quantization) of a polynomial algebra
  is equivalent to a deformation (deformation quantization) whose
  cochains are differential operators.
\end{thm}
\begin{proof}
  Suppose that $\ms=\m+\h\m_1+\dots$ is the given deformation, and
  that $\m_1\mcom\m_n$ are given by differential operators.  Then we
  show that $\m_{n+1}$ can be replaced by a differential operator
  yielding an equivalent deformation.  By Theorem \ref{diffcob}, we
  can express $\delta(\m_{n+1})=\delta(C)$, for some differential
  operator $C$, so that $\delta(\m_{n+1}-C)=0$. Therefore we can
  express $\m_{n+1}-C=A+S$ where $A$ is an antisymmetric cocycle and
  $S$ is a symmetric cocycle. Let $\m'_{n+1}=C+A$.  $A$ is a
  differential operator because it is a biderivation, hence
  $\m_{n+1}'$ is a differential operator.  $\m_{n+1}$ and $\m'_{n+1}$
  differ by $S$, which is a coboundary, since it is a symmetric
  cocycle, so we can replace $\ms$ by an equivalent deformation whose
  first $n+1$ terms are given by differential operators. If $\ms$ is a
  deformation quantization then (\ie, if it is alternating) then the
  new deformation will also be a deformation quantization because the
  differential operator $C$ for which $\delta\pi_{n+1}=\delta C$ has
  the same parity as $\pi_{n+1}$ in view of Theorems \ref{3skew} and
  \ref{3sym}.
\end{proof}
\section{Construction of the cochains $\m_2$ and $\m_3$}
In this section we start {}from a polynomial Poisson algebra $(\A,\Pb)$
over a field $\k$ of characteristic zero and construct an explicit
third order deformation $\m+\h\m_1+\h^2\m_2+h^3\m_3$ of the
multiplication $\m$ in $\A$. In principle the theorems of sections
\ref{hoch} and \ref{hochdiff} allow one to construct a third order
deformation. However, even in the case in which we are given a
concrete example of $\pi_1$ it is difficult to determine $\pi_2$ and
$\pi_3$ explicitly {}from these theorems. Therefore we will use a
different method for constructing $\pi_2$ and $\pi_3$. It should be
remarked that the two constructions do not give the same $\m_2$ and
$\m_3$ terms.

Throughout this section a basis $\{x_i\}_{i\in\I}$ will be fixed and
we use $\X ij$ as a convenient notation for the Poisson bracket
$\{x_i,x_j\}$ and we use superscripts to denote partial derivatives,
as in the previous section. Without loss of generality we pick
$\m_1={\frac 12}\Pb$.  If we use the summation convention then $\pi_1$
can be written as ${\frac12}\X ij\p^i\t\p^j$, the Jacobi identity for
$\pi_1$ reads
\begin{equation}\label{Jacobi_ori}
  \X ij^l\X kl+\X jk^l\X il+\X ki^l\X jl=0
\end{equation}
(for any $i,\,j,\,k\in \I$), the derivative of the Jacobi identity is
written as
\begin{equation}\label{Jacobi_der}
  \X ij^{lm}\X kl+\X jk^{lm}\X il+\X ki^{lm}\X jl
  +\X ij^l\X kl^m+\X jk^l\X il^m+\X ki^l\X jl^m=0,
\end{equation}
(for any $i,\,j,\,k,\,m\in \I$) and there are similar expressions for
higher derivatives. We first give a formula for $\m_2$ and show that
it solves the second equation in (\ref{cohom}).
\begin{prp}\label{order2}
  Given an infinitesimal deformation $\m+\h\m_1$ of $\A$ where
  $\m_1={\frac12}\X ij\p^i\t\p^j$ is antisymmetric and satisfies the
  Jacobi identity, let $\m_2$ be the following symmetric cochain
  \begin{equation}\label{m2}
    \m_2= {\frac1{12}}\X ij^l\X lk\(\p^i\t\p^{jk}+\p^{jk}\t\p^i\)
            +{\frac18}\X ij\X kl\p^{ik}\t\p^{jl}.
  \end{equation}
  Then $\m+\h\m_1+\h^2\m_2$ is a second order deformation of $\A$.
\end{prp}
\begin{proof}
  Use \refeq{GHbracket} to compute the right hand side of
  \begin{equation}\label{eq2}
    \d\pi_2 ={\tfrac12}[\pi_1,\pi_1]
  \end{equation}
  and use the Jacobi identity (\ref{Jacobi_ori}) to find
  \begin{equation}\label{m1m1}
    \tfrac12[\m_1,\m_1]=\tfrac14\X ik^l\X lj\tt ijk+
    \tfrac14\X ij\X kl\(\tt{ik}lj-\tt ik{jl}\).
  \end{equation}
  The third order part of $\d\m_2$ (with $\m_2$ given by (\ref{m2}))
  is computed using (\ref{delta}) to be given by
  \begin{equation*}
    {\frac1{12}}\X ij^l\X lk\(\tt i j k+\tt i k j -\tt j k i -\tt k j i\).
  \end{equation*}
  Since $i,\,j$ and $k$ are just summation indices this can be
  rewritten as $$ {\frac1{12}}\(\X ij^l\X lk+2\X ik^l\X lj- \X kj^l\X
  li\)\tt ijk.  $$ Using the Jacobi identity (\ref{Jacobi_ori}) this
  reduces to a single term $$ {\frac14}\X ik^l\X lj\tt ijk, $$ which
  is the third order term of ${\frac12}[\m_1,\m_1]$. For the fourth
  order term one makes a similar computation (but the Jacobi identity
  is not used).
\end{proof}
One concludes {}from these computations that it is not obvious how to
guess a cochain whose coboundary is given; compare carefully
(\ref{m1m1}) and (\ref{m2}).

Our next task is to find an explicit solution for the the third
equation in (\ref{cohom}), namely the equation
$\d\pi_3=[\pi_1,\pi_2]$. The computation of the right hand side is
long but straightforward. Writing it as a coboundary is non-trivial
and we will concentrate on this aspect.  Clearly every term of the
right hand side is a differential operator of order 3, 4, 5 or 6.  We
will denote the $i$-th order part of a bidifferential operator by a
subscript $(i)$.  We start with the highest order, which is the
easiest.
\begin{lma}\label{order6}
  The sixth order part of $\pp$ is the coboundary of an antisymmetric
  2-cochain,
  \begin{equation}\label{order6form}
    \pp_{(6)}=\d\(\frac1{48}X_{ij}X_{kl}X_{mn}\p^{ikm}\t\p^{jln}\).
  \end{equation}
\end{lma}
\begin{proof}
  The sixth order terms in $\pp$ are the ones for which none of the
  coefficients in $\pi_1$ or $\pi_2$ are differentiated. There are
  twelve terms, they come {}from the bracket of $\m_1$ and the fourth
  order term of $\m_2$ only, and eight of them cancel in pairs,
  leaving the following expression for $\pp_{(6)}$.
  \begin{multline*}
    {\frac1{16}}X_{mn}X_{ij}X_{kl}
          (\tt{ikm}{jl}n+\tt n{jl}{ikm}+\tt{ikm}n{jl}+\tt{jl}n{ikm}).
  \end{multline*}
  To compute $\d(X_{ij}X_{kl}X_{mn}\p^{ikm}\t\p^{jln}),$ use
  (\ref{delta}) and find twelve terms which come in equal triples due
  to the order three symmetry $(i,j)\to(k,l)\to(m,n)$.
  Formula~(\ref{order6form}) follows.
\end{proof}
Note that the computation did not involve the Jacobi identity.  In the
symplectic case this is the only term which survives.  Next, we
consider the terms of order 5.
\begin{lma}\label{order5}
  The fifth order term $\pp_{(5)}$ is also the coboundary of an
  antisymmetric 2-cochain, given by
  \begin{equation*}
    \pp_{(5)}={\frac1{24}}\d\(X_{ij}^kX_{kl}X_{mn}
    \(\p^{jlm}\t\p^{in}-\p^{in}\t\p^{jlm}\)\).
  \end{equation*}
\end{lma}
\begin{proof}
  The bracket $[\pi_1,\pi_2]_{(5)}$ has a lot of terms, they are of
  types \type113, \type131, \type311, \type122, \type212 and \type221.
  The terms of type $\type131$ cancel and in the other ones there is
  some simplification. Since $[\pi_1,\pi_2]$ is flip symmetric and the
  coboundary of any antisymmetric 2-cochain is flip symmetric as well,
  we only need to consider the terms of type \type311, \type122 and
  \type212. We give the result below, omitting a global factor
  $1/24$. Note the non-triviality of the coefficients.
  \begin{align*}\label{types}
    \type311&:\!\(\X nm\X ij^l\X kl+\X im\X nj^l\X kl\)\tt {jkm}in,\\
    \type122&:\!\(\X ij^k\X kl\X mn\!+\!2\X ln^k\X km\X ij\!+\!\X nm^k\X kl\X ij
    +3\X ml^k\X kn\X ij\)\!\tt n{jl}{mi},\\
    \type212&:\!\(\X ij^k\X kl\X mn\!-\!\X ij^k\X km\X ln\!+\!3\X ij^k\X kn\X ml\)
    \tt{im}n{jl}.
  \end{align*}
  It is surprising that all these terms integrate to a single term,
  \ie, as a whole they can be written as
  \begin{equation}\label{int5}
    \d\(\X ij^k\X kl\X mn\(\p^{jlm}\t\p^{in}-\p^{in}\t\p^{jlm}\)\).
  \end{equation}
  Before checking this, note that~(\ref{int5}) produces indeed
  precisely terms of the appropriate types. Clearly the \type311 part
  of~(\ref{int5}) is given by
  \begin{equation*}
    \X ij^k\X kl\X mn\(\tt{jlm}in+\tt{jlm}ni\)
  \end{equation*}
  and is easily rewritten in the form of type \type311. Type \type212
  involves the Jacobi identity. The \type212 part of (\ref{int5}) is
  given by
  \begin{align*}
    -\X ij^k\X kl\X mn
    &(\tt{jl}m{in}+\tt{jm}l{in}+\tt{lm}j{in}\cr
    &+\tt{in}j{lm}+\tt{in}l{jm}+\tt{in}m{jl}),
  \end{align*}
  which is easily rewritten as
  \begin{align*}
    (\X ij^k\X km\X nl&+\X ij^k\X kl\X mn+2\X ij^k\X kn\X ml\\
                     &+\X nj^k\X ki\X ml+\X in^k\X kj\X ml)\tt{im}n{jl}.
  \end{align*}
  Now use the Jacobi identity (\ref{Jacobi_ori}) on the last two terms
  to obtain the term of type \type212. Finally, the \type122 part
  of~(\ref{int5}) is given by $$ -\X ij^k\X kl\X mn\(\tt m{jl}{in}+\tt
  j{lm}{in}+\tt l{jm}{in}\).  $$ When this is rewritten as $$\(\X
  ij^k\X kl\X mn+\X mn^k\X kl\X ij+\X ml^k\X kn\X ij\) \tt n{jl}{mi}
  $$ then the first term matches with the first term of type \type122
  and the other two match up with the three remaining terms of type
  \type122.
\end{proof}
For the fifth order term we used the Jacobi identity.  For the fourth
order term we will also use the derivative of the Jacobi identity
(\ref{Jacobi_der}).
\begin{lma}\label{order4}
  The fourth order term $\pp_{(4)}$ is the coboundary of an
  antisymmetric 2-cochain,
  \begin{align*}
    \pp_{(4)}=
    {\frac1{48}}\d(&\X lm^k\X jn^l\X ki(\p^{mn}\t\p^{ij}-\p^{ij}\t\p^{mn})\\
                +&\X mn^{kl}\X lj\X ki(\p^m\t\p^{nij}-\p^{nij}\t\p^m)).
  \end{align*}
\end{lma}
\begin{proof}
  As in the previous case we give the terms in $\pp_{(4)}$ by
  type. There are just three types, to wit, \type112, \type 121 and
  \type211. By flip symmetry we only need to consider the terms of
  type \type 112 and \type121. They have the following form (we omit
  the global constant 1/48).
  \begin{align*}
    \type112&:\X ki\(4\X ml^k\X jn^l+2\X lj^k\X mn^l+
              2\X nj^{kl}\X lm+3\X mn^{kl}\X lj\)\tt m n{ij},\\
    \type121&:2\(\X mn^{kl}\X lj\X ki-\X ij^{kl}\X km\X ln\)\tt m{ni}j.
  \end{align*}
  We already simplified these formulas by using the Jacobi identity
  (for \type121 we used it twice). The verification for type \type121
  is straightforward: the six terms of type \type121 in
  \begin{equation*}
    \d\(\X mn^{kl}\X ki\X lj\p^m\t\p^{nij} +\X
    ij^{kl}\X km\X ln\p^{mni}\t\p^j\)
  \end{equation*}
   come in pairs and reduce to \type121 above. The terms of type
   \type112 in $$\d\(2\X lm^k\X jn^l\X ki\p^{mn}\t\p^{ij} +\X
   mn^{kl}\X ki\X lj\p^m\t\p^{nij}\) $$ are given by $$\X ki\(-2\X
   ln^k\X jm^l-2\X lm^k\X jn^l +\X mn^{kl}\X lj+2\X mj^{kl}\X ln\)\tt
   m n{ij} $$ which reduces to $$\X ki\(4\X ml^k\X jn^l+2\X lj^k\X
   mn^l+ 2\X nj^{kl}\X lm+3\X mn^{kl}\X lj\)\tt m n{ij} $$ by using
   the derivative of the Jacobi identity.
\end{proof}
Finally we consider the term of order 3. The proof does not involve
the Jacobi identity and is left to the reader.
\begin{lma}\label{order3}
  The third order term $\pp_{(3)}$ is also the coboundary of an
  antisymmetric 2-cochain,
  \begin{equation*}
    \pp_{(3)}={\frac 1{24}}\d(\X ij\X kl^i\X
           mn^{jk}(\p^n\t\p^{lm}-\p^{lm}\t\p^n)).
  \end{equation*}
\end{lma}
Our previous results lead to the following theorem.
\begin{thm}\label{constthm}
  Let $(\A,\Pb)$ be a polynomial Poisson algebra with basis
  $\{x_i\}_{i\in\I}$ and denote $\m_1=\X ij\p^i\t\p^j$, where $\X
  ij=\{x_i,x_j\}$. Then the following formula gives a third order
  deformation $\m+\h\m_1+\h^2\m_2+\h^3\m_3$ of $\A$,
  \begin{align}\label{ourformula}
    \ms=\m&+{\frac{\h}2}\X ij\p^i\t\p^j+{\frac{\h^2}{24}}
        \left[2\X ij^l\X lk\(\p^i\t\p^{jk}+\p^{jk}\t\p^i\)
          +3\X ij\X kl\p^{ik}\t\p^{jl}\right]\notag\\
      &+{\frac{\h^3}{48}}[2\X ij\X kl^i\X mn^{jk}(\p^n\t\p^{lm}-
         \p^{lm}\t\p^n)+X_{ij}X_{kl}X_{mn}\p^{ikm}\t\p^{jln}\\
      &\qquad+\X lm^k\X jn^l\X ki(\p^{mn}\t\p^{ij}-\p^{ij}\t\p^{mn})\notag\\
      &\qquad+\X mn^{kl}\X lj\X ki(\p^m\t\p^{nij}-\p^{nij}\t\p^m)\notag\\
      &\qquad+2X_{ij}^kX_{kl}X_{mn}\(\p^{jlm}\t\p^{in}-\p^{in}\t\p^{jlm}\)].
               \notag
  \end{align}
  Up to equivalence every extension of $\m+\h\m_1$ is of the form
  \begin{equation}\label{gen_def}
    \m+\h\m_1+\h^2(\m_2+\ph_2)+\h^3(\m_3+\ph_3+\psi_3)
  \end{equation}
  with $\ph_2$ and $\ph_3$ antisymmetric biderivations and
  $\psi_3$ a symmetric 2-cochain satisfying
  $\p\psi_3=[\m_1,\ph_2]$. Conversely, for such $\ph_2,\,\ph_3$ and
  $\psi_3$ (\ref{gen_def}) is always a third order deformation.
\end{thm}
\begin{proof}
  We proved already that~(\ref{ourformula}) is a third order
  deformation. Suppose now that $\m+\h\m_1+\h^2\m_2'+\h^3\m_3'$ is
  another deformation which extends the same infinitesimal
  deformation. Then $\ph_2=\m_2-\m_2'$ is a cocycle which can be
  assumed to be an antisymmetric biderivation (Proposition
  \ref{isom}).  Since $\d\m_3'=[\m_1,\m_2+\ph_2]$ is a coboundary its
  flip antisymmetric part $[\m_1,\ph_2]$ satisfies the Jacobi identity
  (Theorem \ref{3sym}); we let $\psi_3$ be any symmetric cochain whose
  coboundary is $[\m_1,\ph_2]$. Then $\m_3'-\psi_3$ must differ {}from
  $\m_3$ by a cocycle $\ph_3$ which we may assume, again without loss
  of generality, to be an antisymmetric biderivation.
\end{proof}
\section{The obstruction to a fourth order deformation}
In this section we want to investigate the fourth order term of the
explicit deformation which is given by (\ref{ourformula}). For a given
polynomial Poisson algebra $(\A,\Pb)$ over a field $\k$ of
characteristic zero we will denote the latter deformation by
$\ms=\m+\h\m_1 +\h^2\m_2+\h^3\m_3$; as before
$\m_1={\frac12}\Pb={\frac12}\X ij\p^i\t\p^j$.
\begin{thm}\label{obsthm}
  The deformation (\ref{ourformula}) extends to a fourth, hence fifth,
  order deformation if and only if the following, non-trivial,
  condition is satisfied for any $\a<\b<\c\in\I$: %
  \begin{align}\label{obs}
    &2\X ij\X kl^i(\X \a\b^{km}\X \c m^{jl}+\X
    \b\c^{km}\X \a m^{jl}+\X \c\a^{km}\X \b m^{jl})\\ &+\X ij\X kl(
    \X\a\b^{ikm}\X\c m^{jl}+\X\b\c^{ikm}\X\a m^{jl}+\X\c\a^{ikm}\X \b m^{jl})=0.
    \notag
  \end{align}
\end{thm}
\begin{proof}
  The deformation (\ref{ourformula}) extends to a fourth order
  deformation if and only if $[\pi_1,\pi_3]+{\frac 12}[\pi_2,\pi_2]$
  is a coboundary. Since this cocycle is flip antisymmetric this is
  equivalent to $J([\m_1,\m_3]+{\frac12}[\m_2,\m_2])=0$. The fact that
  $\m_2$ is symmetric implies at once that $J([\m_2,\m_2])=0$. As for
  the terms in $J([\m_1,\m_3])$, they have orders ranging {}from 3 to 8
  only. We claim that the terms of order at least four all vanish,
  sketching the computation in the least trivial case when the order
  equals four. A direct application of (\ref{GHbracket}) gives the
  following expression for the coefficient of
  $\p^{\a\ab}\t\p^\b\t\p^\c$ in $J([\m_1,\m_3])$ (some indices have
  been relabelled for later convenience and a global constant 1/48 has
  been omitted; note also that $\a$ and $\ab$ can be freely
  interchanged):
  \begin{align*}
    &2\X ji\X l\ab^i(\X\a k\X\b\c^{jkl}+\X\b k\X\c\a^{jkl}+\X\c k\X\a\b^{jkl})\\
    &+2\X ji\X l\ab^i(\X\a k^{lj}\X\b\c^{k}+\X\b k^{lj}\X\c\a^{k}+\
    \X\c k^{lj}\X\a\b^{k})\\
    &+2\X jl\left(\X\ab\b^{ij}(\X\c k\X i\a^{kl}+\X\a k\X\c i^{kl})
    -\X\c\a^{kl}(\X\b i\X k\ab^{ij}+\X\ab i\X\b k^{ij})\right)\\
    &+2\X k\ab^i\X li^j(\X\b j\X\c\a^{kl}+\X\c j\X\a\b^{kl})
    -\X k\ab^i\X\b\c^{kl}(\X\a j\X il^j+\X lj\X i\a^j)\\
    &-2\X li\X jk^i(\X\ab\b^j\X\c\a^{kl}-\X\c\a^j\X\ab\b^{kl}).
  \end{align*}
  We now use the second derivative of the Jacobi identity, i.e., we
  use the formula (valid for any indices $\a,\b,\c,j$ and $l$),
  \begin{equation*}
    (\X\a k\X\c\b^{k}+\X\b k\X\a\c^{k}+\X\c k\X\b\a^{k})^{jl}=0
  \end{equation*}
  to rewrite the first two lines (giving the first line below) and we
  use twice a derivative of the Jacobi identity to rewrite the third
  line (giving lines two and three below); the fourth line is
  simplified by a direct application of the Jacobi identity,
  \begin{align*}
    &2(\X ij\X l\ab^i+\X il\X j\ab^i)
    (\X\a k^j\X\b\c^{kl}+\X\b k^j\X\c\a^{kl}+\X\c k^j\X\a\b^{kl})\\
    &+2\X lj\X\a\b^{ij}(\X\c k^l\X i\ab^k+\X\ab k^l\X\c i^k+\X ik^l\X\ab\c^k)\\
    &-2\X lj\X\c\a^{kl}(\X\b i^j\X k\ab^i+\X\ab i^j\X\b k^i+\X ki^j\X\ab\b^i)\\
    &+2\X k\ab^i\X li^j
    (\X\a j\X\b\c^{kl}+\X\b j\X\c\a^{kl}+\X\c j\X\a\b^{kl})\\
    &-2\X li\X jk^i(\X\ab\b^j\X\c\a^{kl}-\X\c\ab^j\X\a\b^{kl}).
  \end{align*}
  Most terms in this expression cancel out in pairs, leaving
  \begin{align*}
    &2(\X jl\X\c i^j+\X ji\X l\c^j +\X j\c\X il^j)\X k\ab^i\X\a\b^{kl}\\
    &\quad+2(\X jl\X\a i^j+\X j\a \X il^j)\X k\ab^i\X\b\c^{kl}\\
    &\quad+2(\X jl\X\b i^j+\X ji\X l\b^j +\X j\b\X il^j)\X k\ab^i\X\c\a^{kl}
  \end{align*}
  which is zero, by a single application of the Jacobi identity on
  every line. It follows that the only non-zero terms in
  $J([\m_1,\m_3])$ are terms of type $(1,1,1)$. Using
  (\ref{GHbracket}) we find that the coefficient of
  $\p^\a\t\p^\b\t\p^\c$ in $J([\m_1,\m_3])$ is given, (up to a global
  constant $-1/48$) by the left hand side of (\ref{obs}); since this
  expression is antisymmetric in $\a,\,\b,\,\c$ it will hold in
  general when it holds for $\a<\b<\c\in\I$. We will see later an
  example for which (\ref{obs}) is non-zero, showing that our
  deformation (\ref{ourformula}) in general does not extend to a
  fourth order deformation. However, if (\ref{obs}) vanishes then
  $\m_4$ can be chosen to be symmetric, which implies the existence of
  $\m_5$ upon using Theorem \ref{3skew}.
\end{proof}
We will show in Section \ref{remobs} how to overcome this
obstruction. Before doing this we will examine the quantized
enveloping algebra as a natural candidate for a deformation
quantization. As it turns out the same obstruction found in Theorem
\ref{obsthm} will arise. This surprising fact is a consequence of the
non-trivial fact that the third order deformation quantization given
by \reff{ourformula} coincides with the third order deformation
quantization given by the quantized enveloping algebra.
\section{The quantized universal enveloping algebra}
In this section we will show that the third order deformation which we
constructed for any polynomial Poisson algebra comes {}from a
``quantized'' enveloping algebra.  The fact that an enveloping algebra
appears here is not surprising.  The symmetric algebra of a Lie
algebra is a polynomial Poisson algebra in a natural way and it is
well known that the quantized universal enveloping algebra of a Lie
algebra is a deformation quantization of this Poisson algebra
(see~\cite{bffls}, \cite{ber}).

In order to describe the enveloping algebra of a polynomial Poisson
algebra we will view polynomial algebras as symmetric algebras over a
vector space.  Let $V$ be a (possibly infinite-dimensional) vector
space over a field $\k$ of characteristic zero. For simplicity of notation
we will denote elements in
$V$ by lowercase roman letters.
For any positive integer $n$ we let $V^{n}=V\t V\t\dots\t V$ ($n$
copies) and $V^{0}=\k$.  The tensor algebra over $V$ is the
$\Z$-graded associative algebra (with unit) defined by
\begin{equation*}
  \TV=\bigoplus\limits_{n=0}^\infty V^{ n}.
\end{equation*}
The symmetric algebra $\SV$ is the quotient $\SV=\TV/I$, where $I$ is
the homogeneous ideal in $\TV$ generated by elements of the form
$x\tns y-y\tns x$. The symmetric algebra is isomorphic to the
polynomial algebra $\k[x_j]_{j\in \I}$ where $\{x_j\}_{j\in \I}$ is
any basis for $V$. (Of course, any polynomial algebra can be represented in this
form.)  In particular, we will use juxtaposition to denote the
product in $\SV$, just as we did for a polynomial algebra.

Any antisymmetric map $V\t V\to \SV$ extends to a unique antisymmetric
biderivation on $\SV$. When this biderivation satisfies the Jacobi identity
then $(\SV,\Pb)$ becomes a polynomial Poisson algebra, and every
polynomial Poisson algebra arises in this fashion.  The quotient map
$\mu:\TV\to\SV$ has a $\k$-linear right inverse $\sg:\SV\to\TV$ which
is defined by
\begin{equation*}
  \sg\left(\prod_{i=1}^n a_i\right)={\frac 1 {n!}}\sum_{p\in S_n}
        a_{p(1)}\t a_{p(2)}\t\cdots\t a_{p(n)},
\end{equation*}
where $S_n$ is the symmetric group on $n$ elements. We call $\sg$ the {\em
symmetrization map.} Note that $\mu$ is an algebra homomorphism but
the symmetrization map $\sg$ is not.
Let $\TVh$ ($\SVh$) be the formal power series with coefficients in
$\TV$ ($\SV$). Then $\TVh$ and $\SVh$ are naturally $\kh$-algebras,
$\mu$ extends to an $\kh$-algebra homomorphism $\mu:\TVh\ra\SVh$, and
$\sg$ extends to a $\kh$-linear map $\sg:\SVh\ra\TVh$.  Now we
introduce a natural candidate for a deformation quantization of a
polynomial Poisson algebra $(\SV,\Pb)$.
\begin{dfn}
  Let $J^\h$ denote the two-sided ideal of $\TVh$, generated by all
  elements
\begin{equation}\label{rel}
  x\t y-y\t x-\h\sg\{x,y\}\qquad (x,\,y\in V).
\end{equation}
The {\em quantized universal enveloping algebra} of $(\SV,\Pb)$
is given by
\begin{equation}\label{env}
  \Uh=\TVh/J^\h.
\end{equation}
The induced product on $\Uh$ is denoted by $\odot$ and the quotient map by
\begin{equation*}
  \rho:\TVh\to\Uh.
\end{equation*}
\end{dfn}
Thus, we have associated to a polynomial Poisson algebra $(\SV,\Pb)$ a
new (non-commutative) associative algebra $(\Uh,\odot)$ and they are
linked by the $\kh$-linear map (not a homomorphism!)
\begin{equation*}
  \tau:\SVh\to\Uh
\end{equation*}
given by $\tau=\rho\circ\sg$. The maps $\tau,\,\rho$ and $\sg$ induce
maps $\tau_n,\,\rho_n$ and $\sg_n$ on the quotient spaces $\TVn$,
$\SVn$ and $\Un$ obtained by dividing out by the ideal
$(h^{n+1})$. We also use the notation $\Jn$ for
$J^\h/(h^{n+1})$, so that $\Un=\TVn/\Jn$.
We will see that in some important cases the map
$\tau$ is a bijection, but that in general $\tau_n$ is only injective
for $n\leq3$.   If $\tau$ is injective up to some order, the
enveloping algebra provides a deformation quantization of $(\SV,\Pb)$
of the same
order, as given by the following theorem.
\begin{thm}
  If $\tau:\SVh\to\Uh$ (resp.\ $\tau_n$) is injective
  then the unique product $\star$ on $\SVh$ which makes $\tau$ (resp.\
  $\tau_n$) into a homomorphism is a deformation quantization (resp.\
  of order $n$) of the Poisson algebra $(\SV,\Pb)$.
\end{thm}
\begin{proof}
$\tau$ is always surjective: simply note that $\Uh_1$ is
canonically isomorphic to $\SV$,
so that for any $q\in\Uh$ there exists a $p\in\SV$
such that $\tau(p)=q\mod h$. Then $\tau(p)-q=hq_1$, for some
$q_1\in \Uh$. Continuing this process, we obtain a sequence of
polynomials $p_i$ such that $\tau(p+hp_1\mplus h^kp_k)-q=h^kq_k$
for some $q_k\in \Uh$. Then $\tau(p+hp_1+\dots)=q$. It follows
that $\tau_n$ is also surjective.

 If $\tau_n$ is injective then the
associative product which is induced by $\tau_n$ is given for
$p,q\in\SV$ by
\begin{equation*}
  p\star q=\tau_n^{-1}(\tau_n(p)\odot\tau_n(q)).
\end{equation*}
We show that it defines a deformation of $(\SV,\Pb)$ and that it is
alternating. It is easy to see that
\begin{equation*}
  \tau_n(p)\odot\tau_n(q)=\tau_n(pq)\mod h
\end{equation*}
so that $p\star q=pq\mod h$; the associativity of $\star$ on
$\SVn$ implies that $p\star q=pq+\h\m_1(p,q)\mod h^2$ for
some cocycle $\m_1$.  If we can show that $\m_1$ is antisymmetric then
it is a biderivation and the fact that $\m_1={\frac12}\Pb$ follows
{}from the following check for elements $x,y\in V$,
\begin{equation*}
  \h\m_1(x,y)={\frac12}(x\star y-y\star x)={\frac h2}\{x,y\}\mod \h^2.
\end{equation*}
Now we show that $\star$ is alternating (up to order $n$), which
proves in particular that $\m_1$ is antisymmetric.  Let $T$ be the
anti-involution on $\TVh$ induced by the map which reverses the order
of elements in a tensor product, and let $t$ be the involution of $\kh$
which is given by the map $h\mapsto -h$. Then $t$ determines
involutions of $\SVh$ and $\TVh$, which we will also denote by
$t$. Let $\iota=T\circ t=t\circ T$, so $\iota$ is an anti-involution
of $\TVh$. Note that $T\circ\sg=\sg$.  Thus $\iota(x\tns y -y\tns x
-\h\{x,y\})=y\tns x-x\tns y-\h\{y,x\}$, so $\iota$ maps the ideal
$J^\h$ to itself inducing an anti-involution $\imath$.  We also have
the relations $\imath\circ\rho_n=\rho_n\circ\imath$ and $\tau_n\circ
t=\iota\circ \tau_n$.  Now $\star$ is alternating precisely when
$t(p\star q)=q\star p$ for all $p$, $q$ in $\SV$.  But note that
\begin{align*}
  \tau_n(t(p\star q))&= \iota(\tau_n(p\star q))
                      = \iota(\tau_n(p)\odot\tau_n(q))
                      = \iota(\rho_n(\sg_n(p))\odot \rho_n(\sg_n(q)))\\
                     &= \iota(\rho_n(\sg_n(p)\tns\sg_n(q)))
                      = \rho_n(\iota(\sg_n(p)\tns\sg_n(q)))
                      = \rho_n(\sg_n(q)\tns\sg_n(p))
\end{align*}
and similarly $\tau_n(q\star p)=\rho_n(\sg_n(q)\tns\sg_n(p))$.  Since
$\tau_n$ is an isomorphism, the conclusion follows.
\end{proof}
The theorem shows that the injectivity of $\tau_n$ is crucial. We show
in the next theorem how injectivity of $\tau_n$ can be rephrased as an
identity in $\Un$.  Define an antisymmetric map $\Delta:V^{3}\to\Uh$
by
\begin{align*}
  \Delta(x,y,z)&=x\odot\tau\{y,z\}+y\odot\tau\{z,x\}+z\odot\tau\{x,y\}\\
               &-\tau\{y,z\}\odot x-\tau\{z,x\}\odot y-\tau\{x,y\}\odot z
\end{align*}
and call $\Delta=0$ the {\em diamond relation.} For any $n$ there is
an induced map $\Delta_n: V^{\t 3}\to\Un$ and we call $\Delta_n=0$ the
{\em $n$-th diamond relation.}  Note that
for any $x,y,z\in V$,
\begin{equation*}
    \h x\odot\tau\{y,z\}=x\odot y\odot z-x\odot z\odot y,
  \end{equation*}
so that $h\Delta=0$, and similarly $h\Delta_n=0$ for all $n$.
It is precisely the possibility of multiplying a nonzero element in $\Uh$
by $h$ to obtain zero that can cause $\tau$ to fail to be
injective, as we show in the theorem below.
For the proof we need the notion of
ordered elements in the tensor product. Fixing an ordered basis
$\{x_i\}_{i\in\I}$ for $V$ we call an element $\alpha=x_{i_1}\t
x_{i_2}\t\cdots\t x_{i_m}\in\TV$ an {\em ordered} monomial if $i_1\leq
i_2\leq\cdots\leq i_m$, and strictly ordered if the inequalities above
are strict inequalities. Let $\OV$ be the subspace of $\TV$ spanned by the
ordered monomials, $\OVh$ be the induced subspace of $\TVh$,
and $\OVn=\OVh/(h^{n+1})$ be the subspace of ordered elements in $\TVn$.
Also, for an element $\gamma\in \TVn$, denote by $\gamma(0)$
its 0-th order part, so that $\gamma-\gamma(0)\in h\TVn$.
\begin{thm}\label{inj_thm}
  For $n\geq1$ the following four statements are equivalent.
  \begin{enumerate}
    \item $\tau_n$ is injective;
    \item For any $\alpha\in\Un$, $h\alpha=0$ implies $\alpha=0\mod h^n$;
    \item $\star$ satisfies the $n$-th diamond relation $\Delta_n=0;$
    \item The restriction of $\rho_n$ to $\OVn$ is injective.
  \end{enumerate}
Moreover, each of these statements is true for $n=0$.
\end{thm}
\begin{proof}
  Let us first treat the case of $n=0$ because this is used later in
  the proof.  The fact that $\tau_0$ is injective follows immediately
  {}from the fact that the image of $J^\h$ in $\TV$ is the ideal
  $I$, so that $\tau_0$ is essentially the identity map, {}from
  which it also follows that the
  restriction of $\rho_0$ to $\OV_0$ is injective.
  Statements 2) and 3) hold vacuously for $n=0$, so all
  statements are true for $n=0$.

  Let us suppose that $\tau_n$ is injective and let $\alpha\in\Un$ be an
  element such that $h\alpha=0$. Since $\tau_n$ is surjective there exists
  $\beta\in\SVn$ such that $\tau_n(\beta)=\alpha$. Then $\tau_n(h\beta)=
  h\tau_n(\beta) =0$, so that $h\beta=0$ and $\beta\in(h^n)$. Then
  $\alpha=\tau_n(\beta)=0\mod h^n$, which shows that 1) implies 2).

  That 2) implies 3) follows {}from the fact that $h\Delta_n=0$.

  We now show that 4) implies 1), so we assume that the restriction of
  $\rho_n$ to $\OVn$ is injective. We show that $\tau_n$ is
  injective. By induction, we can assume that this theorem is true for
  $n-1$, so that $\tau_{n-1}$ is injective, since $\Delta_{n-1}=0$ if
  $\Delta_{n}=0$.  Therefore, if $\tau_n(\gamma)=0$ for some
  $\gamma\in\SVn$, then since $\tau_{n-1}(\gamma)=0$, we must have
  $\gamma=0\mod h^n$. Thus $\gamma=h^np$ for some $p\in\SV$. But if
  $x_{i_1}\cdots x_{i_k}$ satisfies $i_1\le\cdots\le i_k$, then
  $\tau_n(h^nx_{i_1}\cdots x_{i_k})=h^n\rho_n(x_{i_1}\mtns x_{i_k})$,
  because we can always reorder the terms appearing in a tensor at the
  price of adding $h$ times something. If we express
  $p=\sum_Ia^Ix_{i_1}\cdots x_{i_k}$, where we sum over all increasing
  multi-indices $I=(i_1\mcom i_k)$, and $\beta=\h^n\sum_I
  a^Ix_{i_1}\mtns x_{i_k}$, then $\beta\in\OVn$ and satisfies
  $\rho_n(\beta)=\tau_n(\gamma)=0$, so that $\beta=0$, by injectivity
  of $\rho_n$ on $\OV_n$.  It follows that $p$ must also vanish, and
  thus $\gamma=0$. This shows that 4) implies 1).

  The rest of the proof is devoted to showing that 3) implies
  4). We fix any $n\geq1$ and assume that $\Delta_n=0$. Since the
  kernel of $\rho_n$ restricted to $\OVn$ is $\OVn\cap J_n^\h$, it
  suffices to show that $\OVn\cap\Jn\subseteq \h\Jn$. An arbitrary
  element $\gamma$ of $\ker\rho_n$ is of the form
  $\gamma=\gamma'+\h\gamma''$ where $\gamma', \gamma''\in\Jn$ and
  \begin{equation}\label{c'}
    \gamma'=\sum_{1\le l\le N} \alpha_{l}\tns(x_{i_{l}}\tns
     x_{j_{l}}-x_{j_{l}}\tns x_{i_{l}} -h\sg\{x_{i_{l}},x_{j_{l}}\})\tns\beta_{l}
  \end{equation}
  for some monomials $\alpha_l$, $\beta_l$ in $\TV$, basis elements
  $x_{i_l}$, and $x_{j_l}$ and some positive integer $N$. We need to
  show that if $\gamma$ is ordered then $\gamma'\in\h\Jn$. We first
  show that $\gamma'(0)=0$. Since $\rho_n(\gamma)=0$ also
  $\rho_0(\gamma(0))=0$ which implies that $\gamma(0)=0$ because $\gamma$ and
  hence also $\gamma(0)$ is ordered. Then $\gamma'(0)$ also vanishes because
  $\gamma(0)=\gamma'(0)$.  Now consider a fixed multi-index $I$ and define
  $\gamma'_I$ by \reff{c'} but summing only over those $l$ for which
  the indices in $\alpha_l\t x_{i_l}\t x_{j_l}\t\beta_l$ coincide with the
  ones in $I$ (including multiplicities). Then evidently
  $\gamma_I'(0)=0$. We will show that this implies that $\gamma_I'\in\h\Jn$, {}from
  which it follows that $\gamma'\in\h\Jn$ because $\gamma'=\sum_I\gamma'_I$.

  First we consider the case when $I$ is a strictly ordered monomial,
  in which case we may assume that $I=(1\mcom m)$ for some $m$. We
  denote by $S_m$ the symmetric group and we consider its standard
  presentation with generators $\theta_k,\,k=1,\dots,m-1$, ($\theta_k$
  corresponds to the transposition $(k,k+1)$) and relations
  $\theta_k^2,\, (\theta_l\theta_{l+1})^3$ and $(\theta_i\theta_j)^2$
  for $\vert i-j\vert\geq2$. For $\lambda\in S_m$, let
  $x_\lambda=x_{\lambda(1)}\mtns x_{\lambda(m)}$.  Then we may express
  $\gamma_I'$ as
  \begin{equation}\label{boundary_sym}
    \gamma_I'=\sum_{\lambda\in S_m}\sum_{k=0}^{m-1}a_{\lambda,k}
    \left(x_{\lambda}-x_{\theta_k\lambda}-\h\chi_\l\right)
  \end{equation}
  where $\alpha_{\lambda,k}\in\k$ and $\chi_{\l,k}=x_{\l(1)}\mtns\sg
  \{x_{\l(k)},x_{\l(k+1)}\}\mtns x_{\l(m)}$. Now consider the Cayley
  graph $\Gamma_m$ of the above presentation for $S_m$. The vertices
  of $\Gamma_m$ are given by the elements in $S_m$, with an edge
  connecting two vertices precisely when the permutations defining
  them differ by a transposition.  The oriented edge connecting
  $\lambda$ and $\theta_k\lambda$ is denoted by $e_{\lambda,k}$, so
  that $\p(e_{\lambda,k})=\lambda- \theta_k\lambda$. We define a
  linear map $\Psi$ {}from the group $C^1(\Gamma_m,\k)$ of (oriented)
  1-chains on $\Gamma_m$ to $\TVh$ by letting
  \begin{equation*}
    \Psi(e_{\l,k})=x_\l-x_{\theta_k\l}-h\chi_{\l,k}.
  \end{equation*}
  Notice that $\Psi$ is well-defined because although $e_{\theta_k\l,k}$ is the
  same edge as $e_{\l,k}$ but with the opposite orientation, it gets
  mapped to $-\Psi(e_{\l,k})$. Then obviously
  \begin{equation*}
    \gamma_I'=\Psi\left(\sum_{\lambda\in
    S_m}\sum_{k=0}^{m-1}a_{\lambda,k}e_{\l,k} \right)
  \end{equation*}
  and the fact that $\gamma_I'(0)$ vanishes means that $\sum_{\lambda\in
  S_m}\sum_{k=1}^{m-1}a_{\lambda,k}e_{\lambda,k}$ is a cycle in the
  homology of the Cayley graph. By the universal coefficient theorem,
  every cycle (with coefficients in an arbitrary group) on a graph can
  expressed as a sum of multiples of closed edge paths in the graph;
  moreover, any cycle on the Cayley graph of a presentation is a sum
  of cycles (with integral coefficients) which correspond to the basic
  relations which appear in the presentation. It follows that
  $\sum_{\lambda\in
  S_m}\sum_{k=1}^{m-1}a_{\lambda,k}e_{\lambda,k}=\sum_{l=1}^tb_lr_l$
  where each $r_l$ corresponds to one of the basic relations appearing
  in the presentation and $\beta_l\in\k$.  Therefore we have that
  \begin{equation*}
    \gamma_I'=\sum_{l=1}^tb_l\Psi(r_l),
  \end{equation*}
  and it suffices to show that $\Psi(f)\in h\Jh$ for any cycle $f$ which
  corresponds to a basic relation. First, notice that the cycle $f$
  which corresponds to $\theta_k^2$ is zero because it consists of the
  sum of two copies of an edge with opposite orientation. Second, let
  $i$ and $j$ be such that $\vert i-j\vert>1$ and let $f_{ij}$ be the
  corresponding cycle, $f_{ij}=e_{\l,i}+e_{\theta_i\l,j}
  +e_{\theta_j\theta_i\l_i}+e_{\theta_j\l,j}$. Then
  \begin{equation*}
    \Psi(f_{ij})=-\h(\chi_{\l,i}+\chi_{\theta_i\l,j}+\chi_{\theta_j\theta_i\l,i}
    +\chi_{\theta_j\l,j}).
  \end{equation*}
  Now both $\chi_{\l,i}+\chi_{\theta_j\theta_i\l,i}$ and
  $-\chi_{\theta_i\l,j}-\chi_{\theta_j\l,j}$ are given, up to an
  element of $\Jn$, by
  \begin{equation*}
    x_{\l(1)}\mtns\{x_{\l(i)},x_{\l(i+1)}\}\mtns\{x_{\l(j)},
    x_{\l(j+1)}\}\mtns x_{\l(m)},
  \end{equation*}
  showing that $\Psi(f_{ij})\in\h\Jn$. Finally, let us assume that
  $f_l$ corresponds to the relation $(\theta_l\theta_{l+1})^3$. Then
  \begin{equation*}
    f_l=e_{\l,l}+e_{\theta_{l}\l,l+1}+e_{\theta_{l+1}\theta_{l}\l,l}
        +e_{\theta_l\theta_{l+1}\theta_l\l,l+1}
        +e_{\theta_l\theta_{l+1}\l,l}
        +e_{\theta_{l+1}\l,l+1}
  \end{equation*}
  so that
  \begin{multline*}
    \Psi(f_l)=\h x_{\l(1)}\mtns(x_{\l(l)}\t\sg\{x_{\l(l+1)},x_{\l(l+2)}\}\\
           -\sg\{x_{\l(l+1)},x_{\l(k+2)}\}\t x_{\l(k)}+\cycl)\mtns x_{\l(m)}.
  \end{multline*}
  Since $\Delta_n=0$ the term between parentheses lies in $J^h_{n-1}$.
  But now note that if $\alpha\in J^h_{n-1}$, then $\alpha=\beta +h^n\gamma$
  for some $\beta\in \Jn$, so that $h\alpha\in \h\Jn$. Thus we
  can conclude that $\Psi(f_l)\in\h\Jn$.

  This completes the proof that 3) implies 4) in case $I$ is strictly
  ordered. If $I=(i_1,\dots,i_m)$ is merely ordered then the proof can
  repeated verbatim after replacing $S_m$ with a quotient group, whose
  presentation is obtained {}from the above standard presentation of
  $S_m$ by adding the relations $\theta_k$ for any $k$ for which
  $i_k=i_{k+1}$. The corresponding Cayley graph is obtained {}from the
  one for $S_m$ by collapsing the edges which correspond to those
  $\theta_k$.
\end{proof}
The above theorem gives us an analytic criterion to check injectivity
at some order. When we assume that injectivity at order $n-1$ has been
checked then we may think of the $n$-th diamond relation as being a
relation in $\SVn$. Since this is the way in which we will use the
diamond relation below, we formulate this fact in a separate theorem.
\begin{thm}
If $\tau_{n}:\SVn\to\Un$ is injective (hence bijective) then $\tau_{n+1}$
is also injective if and only if the diamond relation
\begin{equation*}
  x_\a\star\{x_\b,x_\c\}-\{x_\b,x_\c\}\star x_\a+\cycl(\a,\b,\c)=0
\end{equation*}
holds for any $\a,\b,\c\in\I$. In this formula $\star$ is the product on $\SVn$
which is induced using $\tau_{n}$.
\end{thm}
In this formulation the theorem will turn out to be very useful.  For
example we note that $p\star q=q\star p\mod\h$ and conclude {}from it
that $\tau_1$ is injective. In order to use the theorem to prove
injectivity of the higher $\tau_i$ we need an explicit formula for the
$\star$-bracket which comes {}from the enveloping algebra. We will
show now that such a formula is given exactly by (\ref{ourformula})
and derive injectivity of $\tau_2$ and $\tau_3$ {}from it.
\par
Given a deformation $(\SVh,\star)$ of $\SV$
there is, besides the enveloping algebra $\Uh$ another
(in general) enveloping algebra which is associated
to it.
\begin{dfn}
Let $(\SVh,\star)$ be a deformation (of finite order or formal) of
$\SV$ and denote the commutator in $(\SVh,\star)$
by $\Lb_\star$. Define $\Jhs$ to be the two-sided ideal of
$\TVh$, generated by all elements of the form
  $$a\t b-b\t a-\sg[a,b]_\star, \qquad(a,b\in V)
  $$
and define the \emph{$\star$-enveloping algebra} $\Uhs$ of
$(\SVh,\star)$ by
  $$ \Uhs=\TVh/\Jhs.
  $$
\end{dfn}
For a given deformation $(\SVh,\star)$ the enveloping  algebras $\Uh$
and $\Uhs$ coincide if and only if
\begin{equation}\label{bracket_exact}
[x,y]_\star=\h\{x,y\}\qquad (x,y\in V).
\end{equation}
We call a deformation which satisfies (\ref{bracket_exact}) {\em bracket-exact.}
In terms of the cocycles $\m_i$ this means that
\begin{equation*}
\m_i(x,y)=0\qquad (x,y\in V,\,i>1).
\end{equation*}
For example, our general formula (\ref{ourformula}) defines a
bracket-exact deformation quantization; adding any non-zero
antisymmetric biderivation to $\m_3$ defines a deformation
quantization which is not bracket-exact.
\par
We now give a property which characterizes $\star$-enveloping algebras; in
the case of bracket-exact deformations it characterizes enveloping algebras,
showing that the $\star$-product which comes {}from the enveloping algebra is
given by (\ref{ourformula}).
\begin{dfn}
Let $(\SVh,\star)$ be a deformation of $\SV$. The $\kh$-linear map,
\begin{equation*}
\sg_\star:\SVh\to\SVh
\end{equation*}
which is defined by
\begin{equation*}
\sg_\star\(\prod_{i=1}^n a_i\)=\frac1{n!}\sum_{p\in S(n)}a_{p(1)}\star a_{p(2)}
  \star\cdots\star a_{p(n)}.
\end{equation*}
is called \emph{$\star$-symmetrization.}  We will say that $\star$ is
{\em $s$-balanced} if $\sg_\star$ is the identity when restricted to
elements of $\SV$ of degree $\leq s$. If $(\SVh,\star)$ is a
deformation (of order $n$) of $\SV$ then we call it a {\em
balanced deformation} if $\star$ is $s$-balanced, where $s$ is the
degree of $\Lb_\star$, i.e., the supremum of the degrees of all
coefficients of $[x,y]_\star$, where $x,\,y$ run over $V$
(this degree may be infinite).
\end{dfn}
Note  that when a deformation is bracket-exact then
the degree of $\Lb_\star$ is the degree of the corresponding Poisson bracket
$\Pb$.
\par
\begin{exa} Any deformation is equivalent to a 2-balanced deformation. Indeed,
such an equivalence is given precisely by  $\sg_\star$, i.e., define an
equivalent product $\circ$ by
\begin{equation*}
  p\circ q=\sg_\star^{-1}(\sg_\star(p)\star\sg_\star(q)).
\end{equation*}
Then
\begin{equation*}
  \sg_\circ(xy)={\frac12}(x\circ y+y\circ
  x)={\frac12}\sg_\star^{-1}(x\star y +y\star x)=xy,
\end{equation*}
for any $x,\,y\in V$, so that $\circ$ is 2-balanced.
\end{exa}
\begin{lma}
  Formula (\ref{ourformula}) gives, for any polynomial Poisson
  algebra, a bracket-exact balanced deformation of order 3.
\end{lma}
\begin{proof}
The proof of balancing is by induction.
Obviously any deformation is 1-balanced, so we assume that the
deformation, given by Formula \reff{ourformula}, is $n$-balanced and
prove that it is $(n+1)$-balanced. To do this, take a monomial $a$ of
degree $n+1$ and write $a=a_1a_2\cdots a_{n+1}$. We denote the
associative product (\ref{ourformula}) on $\SVh_3$ by $\star$ and
the corresponding cochains by $\m_i$.  Using the associativity of
$\star$ one has
\begin{equation*}
  \sum_{\tau\in S_{n+1}}a_{\tau(1)}\star a_{\tau(2)}\star\cdots\star\a_{\tau(n+1)}=
  \sum_{i=1}^{n+1}a_i\star\left(\prod_{j\neq i}^{n+1}a_j\right)
\end{equation*}
so $\star$ is $(n+1)$-balanced when
\begin{equation*}
\sum_{i=1}^{n+1}\m_k\left(a_i,\prod_{j\neq i}a_j\right)=0,
\end{equation*}
for $k=1,\,2,\,3$. The verification is immediate.
\end{proof}
The following theorem gives a precise relation between balanced
deformations and the $\star$-enveloping algebra.
\begin{thm}\label{isom_thm}
If $(\SVh,\star)$ is a balanced deformation of
$\SV$ then the $\kh$-algebra homomorphism
  $$F:(\TVh,\t)\to(\SVh,\star)
  $$
which is induced by the natural inclusion $V\to\SV$ induces an $\kh$-algebra
isomorphism
\begin{equation*}
 f:(\Uhs,\odot)\to(\SVh,\star).
\end{equation*}
When $(\SVh,\star)$ is moreover bracket-exact then $\Uhs=\Uh$
and we have an isomorphism
\begin{equation*}
 f:(\Uh,\odot)\to(\SVh,\star).
\end{equation*}
The corresponding statements for $n$-th order deformations also hold.
\end{thm}
\begin{proof}
We will only prove the first statement. If we denote the canonical map
$\TVh\to \Uhs$ by $\rho_\star$ then it suffices to prove that
$\ker F=\ker\rho_\star$ and that $F$ is surjective.  Let us first show
that $F$ is surjective. If $p\in\SV$ then there exists an element
$\alpha\in\TV$ such that $p=F(\alpha)\mod\h$. Indeed, since $\star$
is a deformation we have for any monomial $\prod_{i=1}^na_i$ that
\begin{equation*}
  \prod_{i=1}^na_i=a_1\star a_2\star\dots\star a_n
                  =F(a_1\t a_2\t\dots\t a_n)\mod\h.
\end{equation*}
More generally, for any $k\in\N$,
since $F$ is $\kh$-linear
  we can find $\alpha_0,\dots,\alpha_k\in\TV$ such that
$p=F(\alpha_0+\alpha_1\h+\cdots+\alpha_k\h^k)\mod\h^{k+1}$. It follows
that $\SV\subset\Im F$, which is sufficient to prove that $F$ is surjective.

Let us show that $\ker\rho_\star=\ker F$. Take $a,\,b\in V$ and compute
\begin{align*}
F(a\t b-b\t a-\sg[a,b]_\star)
                    &=F(a)\star F(b)-F(b)\star F(a)-F\sg[a,b]_\star\\
                    &=a\star b-b\star a-\sg_\star[a,b]_\star\\
                    &=a\star b-b\star a-[a,b]_\star,
\end{align*}
which is zero; we used in the computation that $\sg_\star=F\sg$ and
that $\sg_\star[a,b]_\star=[a,b]_\star$ (because the deformation is
balanced).  This shows that $\ker\rho_\star\subset\ker F$.

To show that $\ker F\subset\ker\rho_\star$ we pick
any $X\in\TVh$ for which $F(X)=0$ and show the existence of
$Y\in\TVh$ such that $\rho_\star(X)=\rho_\star(Y)$ and
whose degree (in $\h$) is larger than the degree of $X$.
This will imply that for any $j\in\N$ the composition
  $$\TVh\,{\buildrel\rho_\star\over\longrightarrow}\,\Uhs\longrightarrow\Uhs/(h^j)
  $$
maps $X$ to 0, hence $\rho_\star(X)=0$. To prove it, let $d$ denote the
degree of $X$, i.e., $X=X_0h^d\mod \h^{d+1}$. Let $\bar X_0$ denote the unique element
in $\Im \sg$ for which
  $$\rho_\star(X_0)=\rho_\star(\bar X_0)\mod \h
  $$
If we write
  $$\bar X_0={\frac c{n!}}\sum_{p\in S_n} a_{p(1)}\t a_{p(2)}\t\cdots\t
  a_{p(n)}
  $$
then
\begin{align*}
F(\bar X_0)
      &={\frac c{n!}}\sum_{p\in S_n} F(a_{p(1)})\star
                           F(a_{p(2)})\star\cdots\star F(a_{p(n)})\\
      &={\frac c{n!}}\sum_{p\in S_n} a_{p(1)}\star a_{p(2)}\star\cdots\star  a_{p(n)}\\
      &=ca_1a_2\dots a_n\mod\h.
\end{align*}
Thus $F(X)=0$ implies that $c=0$ so that $\bar X_0=0$. So there exists
a $Y_0$ such that $\rho_\star(X_0)=\rho_\star(hY_0)\mod\h^2$ and
hence there exists an element $Y\in\TVh$ of the form
$Y=Y_0h^{d+1}\mod\h^{d+2}$ such that $\rho_\star(X)=\rho_\star(Y)$.
\end{proof}
We have seen that Formula (\ref{ourformula}) defines a bracket-exact
balanced deformation (of order three). Theorem \ref{isom_thm} implies
that this deformation comes {}from the enveloping algebra, via the
symmetrization map. This fact has the important consequence that we
can use (\ref{ourformula}) to check injectivity of the maps
$\tau_n$. We used already the first term of our formula; \ie, we have
used $p\star q=pq\mod \h$ to show that $\tau_1$ is injective. Further,
\begin{align*}
&\m_1(x_\a,\{\b,\c\})-\m_1(\{\b,\c\},x_\a)+\cycl(\a,\b,\c)\\
&\qquad
=\{x_\a,\{x_\b,x_\c\}\}+\{x_\b,\{x_\c,x_\a\}\}+\{x_\c,\{x_\a,x_\b\}\}
\end{align*}
which is zero in view of the Jacobi identity. This proves injectivity of $\tau_2$. Also
\begin{equation*}
\m_2(x_\a,\{\b,\c\})-\m_2(\{\b,\c\},x_\a)+\cycl(\a,\b,\c)=0
\end{equation*}
since $\m_2$ is symmetric, hence $\tau_3$ is also injective. The fact
that this step is easy is similar to the fact that the existence of
$\m_3$ is automatic (given the fact that $\m_1$ is antisymmetric and
that $\m_2$ is symmetric). Finally, let us examine the injectivity
of~$\tau_4$.
\begin{align*}
&\m_3(x_\a,\{\b,\c\})-\m_3(\{\b,\c\},x_\a)+\cycl(\a,\b,\c)\\
&\quad={\frac1{24}}(2\X ij\X kl^i\X m\a^{jk}\X\b\c^{lm}+\X \a n^{kl}\X lj\X ki\X\b\c^{nij})
+\cycl(\a,\b,\c)\\
&\quad={\frac1{12}}\X ij\X kl^i(
\X \a\b^{km}\X \c m^{jl}+\X \b\c^{km}\X \a m^{jl}+\X \c\a^{km}\X \b m^{jl})\\
&\quad\qquad
{\frac1{24}}\X ij\X kl(
\X\a\b^{ikm}\X\c m^{jl}+\X\b\c^{ikm}\X\a m^{jl}+\X\c\a^{ikm}\X \b m^{jl}).
\end{align*}
which is identical to the obstruction (\ref{obs}) which we found when
trying to extend the deformation given by (\ref{ourformula}).
We will see in the examples that in general the obstruction is
non-zero, hence $\tau_4$ is not injective and the enveloping algebra
leads in general only to a deformation of order three.
\section{The extension to a fourth order deformation}\label{remobs}
We now come to the existence question of a fourth order deformation
for a polynomial Poisson algebra $(\A,\Pb)$ over a field $\k$ of
characteristic zero. We denote the third order deformation
quantization that we obtained in \reff{ourformula} by
$\ms=\m+\h\m_1+\h^2\m_2+\h^3\m_3$ where $\m_1=\frac12\Pb$.  We have
shown in Theorem \ref{constthm} that we get up to equivalence all
possible third order deformations of $(\A,\Pb)$ by adding any
biderivations $\ph_2$ and $\ph_3$ to $\m_2$ and $\m_3$ and adding any
symmetric cochain $\psi_3$ satisfying $\delta\psi_3=[\m_1,\ph_2])$ to
$\m_3$.  Let us denote such an alternative deformation by
$\ms'=\m+\h\m_1+\h^2\m_2'+\h^3\m_3'$. If $\m_\star'$ extends to a
fourth order deformation by adding a term $\h^4\m_4$ then $\m_4$ is a
solution to
\begin{equation*}
  \d\m_4=[\m_1,\m_3']+{\frac12}[\m_2',\m_2'],
\end{equation*}
and the antisymmetric part of the right hand side must be in the kernel of
$J$, leading to
\begin{equation}\label{obs_gen}
  J\left([\m_1,\m_3]+[\m_1,\ph_3]+{\tfrac12}[\ph_2,\ph_2]\right)=0.
\end{equation}
In view of the following lemma, all terms in the left hand side of
(\ref{obs_gen}) are of type \type111.
\begin{lma}
  If $\ph$ and $\psi$ are two biderivations then $J([\ph,\psi])$ has
  type $(1,1,1)$.
\end{lma}
\begin{proof}
  Let $\ph=Y_{ij}\p^i\t\p^j$ and $\psi=Z_{kl}\p^k\t\p^l$. Then the
  piece of $[\ph,\psi]$ that does not contain terms of type $(1,1,1)$
  is given by
  \begin{align*}
    (Y_{ij}Z_{kl}+Y_{kl}Z_{ij})(\tt{ik}lj-\tt ik{lj}).
  \end{align*}
  Applying the Jacobi map every term appears twice with opposite signs
  hence they all cancel out.
\end{proof}
By computing the terms of type type \type111 in \reff{obs_gen} we find
that the existence of a fourth order deformation for a given
$(\A,\Pb)$ is equivalent to the existence of two antisymmetric
biderivations $\ph_2={\frac1{4}}Y_{ij}\p^i\t\p^j$ and
$\ph_3={\frac1{48}}Z_{ij}\p^i\t \p^j$ such that for any
$\a<\b<\c\in\I$
\begin{align}\label{obssol}
  &\X m\c Z_{\a\b}^m+Z_{m\c}\X\a\b^m+6Y_{m\c}Y_{\a\b}^m-
  \X ij\X kl\X\a\b^{ikm}\X\c m^{jl}-2\X ij\X kl^i\X a\b^{km}\X\c m^{jl}\\
  &\qquad + \hbox{ cycl }(\a,\b,\c)=0.\notag
\end{align}
\begin{lma}
  The 2-cocycles $Y_{\a\b}=0$ and
  \begin{equation}\label{corterm}
    Z_{\a\b}=\frac12\X\a\b^{ik}\X ij^l\X kl^j-\X\a i^{jk}\X\b j^{il}\X kl,\quad(\a,\b\in\I)
  \end{equation}
  solve equation \reff{obssol} hence yield the correction term
  \begin{equation*}
    \ph_3=\frac1{96}(\X mn^{ik}\X ij^l\X kl^j-2\X mi^{jk}\X nj^{il}\X kl)\p^m\t\p^n
  \end{equation*}
  to $\m_3$ in \reff{ourformula} in order for the deformation
  quantization to extend to a fourth order deformation quantization.
\end{lma}
\begin{proof}
  Consider the following four equations, which are all a consequence of the Jacobi identity.
  \begin{align*}
    1/2(\X\a\b^i\X\c i)^{jl}\X jk^m\X lm^k+\hbox{ cycl }(\a,\b,\c)=0,\\
    (\X\a\b^j\X jk+\X\b k^j\X j\a+\X k\a^j\X j\b)^{il}\X\c i^{km}\X lm+
        \hbox{ cycl }(\a,\b,\c)=0,\\
    (\X\c i^j\X jk+\X ik^j\X j\c+\X k\c^j\X ji)^l\X\a\b^{im}\X lm^k+
        \hbox{ cycl }(\a,\b,\c)=0,\\
    (\X\c i^j\X jk+\X ik^j\X j\c+\X k\c^j\X ji)\X \a l^{im}\X \b m^{kl}+
        \hbox{ cycl }(\a,\b,\c)=0.
  \end{align*}
  Expand now $\X m\c Z_{\a\b}^m+Z_{m\c}\X\a\b^m+\hbox{ cycl
  }(\a,\b,\c)=0$, (where $Z_{\a\b}$ is given by \reff{corterm}) and
  add the above four equations. After the smoke clears up you will
  find
  \begin{equation*}
    \X ij\X kl\X\a\b^{ikm}\X\c m^{jl}+2\X ij\X kl^i\X a\b^{km}\X\c m^{jl}
  \end{equation*}
  as needed to solve \reff{obssol}.
\end{proof}
\section{Examples}
In this section we will investigate some general and some more specific
examples.  We use the diamond relations to show that for constant and
linear brackets the quantized enveloping algebra always gives a formal
deformation quantization. For the quadratic case we give a few
examples in which the quantized enveloping algebra gives a fifth order
deformation (at least) and we give an example in which the quantized
enveloping algebra gives a formal deformation quantization. We give in
the cubic case a few examples for which the quantized enveloping
algebra gives a deformation of order three but not of higher order
thereby showing the non-injectivity of $\tau_4$ in general. All these
examples are in $\k^4$ (with coordinates $x_1,\dots,x_4$; $\k$ is a
field of characteristic zero) but they have higher-dimensional
counterparts. We will describe the Poisson structure by a $4\times4$
matrix whose $(i,j)$-th entry is the Poisson bracket $\{x_i,x_j\}$. We
refer to this matrix as the \emph{Poisson matrix}.

The simplest case is the one in which all $\X ij$ are constant (i.e.,
they belong to $\k$). It is well-known that in this case a deformation
quantization always exists. This follows also immediately {}from the
diamond relations: since in this case
\begin{equation*}
  x\odot\tau\{y,z\}-\tau\{y,z\}\odot x=0
\end{equation*}
for any $x,y,z\in V$ we conclude that $\Delta=0$ hence that $\tau$ is
injective.  Alternatively it is immediate to check that the following
explicit formula defines a deformation quantization in this case,
\begin{equation*}
  \m_\star=\m+\sum_{n=1}^\infty{\frac{\h^n}{2^nn!}}\X{k_1}{l_1}\cdots\X{k_n}{l_n}
  \p^{k_1\dots k_n}\t\p^{l_1\dots l_n}.
\end{equation*}
If a linear map $V\t V\to V$ satisfies the Jacobi identity then its
extension to $\SV$ also satisfies the Jacobi identity, hence a Lie
algebra leads in a natural way to a polynomial Poisson algebra. We
call it \emph{linear} because the bracket of any two basis elements is
a linear combination of the basis elements. In this case it is known
that the quantized enveloping algebra defines a formal deformation
quantization. This is checked immediately using the diamond relations:
in this case the fact that $\{y,z\}\in V$ for any $y,z\in V$ implies
that
\begin{equation}\label{diahere}
  x\odot\tau\{y,z\}-\tau\{y,z\}\odot x=\h\{x,\{y,z\}\}
\end{equation}
so that the diamond relation holds in view of the Jacobi identity.
Note also that, as a corollary of Theorem \ref{inj_thm} all
bracket-exact deformations of a linear bracket are isomorphic (to the
one given by the enveloping algebra).

We can also consider brackets which have both linear and constant
terms. Since the constant terms define a central extension of the
linear terms this case is also covered by the linear case and the
quantum enveloping algebra defines a deformation
quantization. Alternatively, it is easy to see that \reff{diahere}
also holds in this case so that again the diamond relation is
satisfied.

A major source of examples of non-linear polynomial Poisson brackets
can be found on page 70 of \cite{van}. Consider $\C^{2d}$ as the
linear space of pairs of polynomials $(u(\l),v(\l))$ with $u(\l)$
monic of degree $d$ and $v(\l)$ of degree less than $d$. If we write
\begin{align*}
  u(\l)&=\l^d+u_1{\l}^{d-1}+\cdots+u_{d-1}{\l}+u_d,\\
  v({\l})&=v_1{\l}^{d-1}+\cdots+v_{d-1}{\l}+v_d,
\end{align*}
then the following formula defines for any polynomial $\ph$ in two
variables a Poisson bracket on $\C^{2d}$,
\begin{equation}\label{mypb}
\begin{split}
\{u({\l}),u_j\}&=\{v({\l}),v_j\}=0,\\
\{u_j,v({\l})\}&=\varphi
    (\l,v(\l))\left[\frac{u(\l)}{\l^{d-j+1}}\right]_+\mod u(\l),
                    \qquad 1\leq j\leq d.
\end{split}
\end{equation}
The subscript $+$ means take the polynomial part and the expression
$p(\l)\mod u(\l)$ means take the remainder obtained by Euclidean division.
Since in these particular examples the Poisson matrix is always of the
form $\left(
\begin{array}{cc}
  0&U\\-U&0\\
\end{array}
\right)$ we will only give the matrix $U$ and the polynomial it derives
{}from. Let us explain shortly how to compute $U$ {}from \reff{mypb} for a
given bracket $\ph$ on $\C^4$.  The coordinates are $u_1,u_2,v_1$ and
$v_2$; also $u(\l)=\l^2+u_1\l+u_2$ and $v(\l)=v_1\l+v_2$. Then the
first row of $U$ consists of the coefficients of $\ph(\l,v(\l))\mod
u(\l)$ (just do Euclidean division) and the second row is given by the
coefficients of $\ph(\l,v(\l))(\l+u_1)\mod u(\l)$.  For example, take
$\ph=x^3$. Then
\begin{equation*}
  U=\left(
  \begin{array}{cc}
    u_1^2-u_2&u_1u_2\\ u_1u_2&u_2^2\\
  \end{array}
  \right).
\end{equation*}
In this case direct substitution in the left hand side of (\ref{obs})
gives zero so that the deformation, as given by (\ref{ourformula}),
extends to a fifth order deformation.  Another quadratic bracket is
found by taking $\ph=y$. Then $U$ is given by
\begin{equation*}
  U=\left(
  \begin{array}{cc}
    v_1&v_2\\ v_2& u_1v_2-u_2v_1\\
  \end{array}
  \right).
\end{equation*}
Again (\ref{obs}) is satisfied.  The same is also true for the sum,
$\ph=x^3+y$, which corresponds to taking the sum of the above $U$
matrices.  Another quadratic example of interest is the quadratic
bracket on $\gl(2)$ (see \cite{lipar}). It has a Poisson matrix
\begin{equation*}
  U=\left(
  \begin{array}{cccc}
  0&x_1x_2&0&x_2x_3\\-x_1x_2&0&0&x_2x_4\\
           0&0&0&0\\-x_2x_3&-x_2x_4&0&0\\
  \end{array}
  \right).
\end{equation*}
(\ref{obs}) is satisfied and the deformation extends to order five.
In the following example of a quadratic bracket the quantized
universal enveloping algebra gives a formal deformation
quantization. If $(a_{ij})$ is a skew-symmetric matrix of size $4$
then $\{x_i,x_j\}=a_{ij}x_ix_j$ defines a quadratic Poisson bracket on
$\C^4$.  In this case the relation
\begin{equation*}
  x_i\odot x_j-x_j\odot x_i=\h\tau\{x_i,x_j\}=\h a_{ij}
            (x_i\odot x_j+x_j\odot x_i)
\end{equation*}
can be rewritten as $x_j\odot x_i=A_{ij}x_i\odot x_j$ where
$A_{ij}=(1-\h a_{ij})/ (1+\h a_{ij}).$ The verification of diamond
relation then reduces to the following computation.
\begin{align*}
  x_i&\odot\{x_j,x_k\}-\{x_j,x_k\}\odot x_i+\cycl(i,j,k)\\
  &=x_i\odot x_j\odot x_k(a_{jk}-a_{ij})+x_i\odot x_k\odot x_j(a_{jk}-a_{ki})+
       x_j\odot x_i\odot x_k(a_{ki}-a_{ij})\\
  &+x_j\odot x_k\odot x_i(a_{ki}-a_{jk})+x_k\odot x_i\odot x_j(a_{ij}-a_{ki})
  +x_k\odot x_j\odot x_i(a_{ij}-a_{jk})\\
  &=x_i\odot x_j\odot x_k((a_{jk}-a_{ij})+(a_{jk}-a_{ki})A_{jk}+(a_{ki}-a_{ij})A_{ij}\\
  &+(a_{ki}-a_{jk})A_{ik}A_{ij}+(a_{ij}-a_{ki})A_{ik}A_{jk}+(a_{ij}-a_{jk})
  A_{ij}A_{ik}A_{jk})\\
  &=0.
\end{align*}
Therefore the quantized enveloping algebra of this quadratic Poisson
bracket gives a formal deformation quantization.

Next we consider a few higher order brackets. As in the quadratic
case, if you take $\ph=x^4$ then
\begin{equation*}
  U=\left(
  \begin{array}{cc}
    -u_1^3+2u_1u_2&u_2^2-u_1^2u_2\\ u_2^2-u_1^2u_2&-u_1u_2^2\\
  \end{array}
  \right).
\end{equation*}
In this case we find again that (\ref{obs}) is satisfied so that the
enveloping algebra leads to a fifth order deformation. However, if you
take $\ph=y^2$ then $U$ is given by
\begin{equation*}
  U=\left(
  \begin{array}{cc}
    2v_1v_2-u_1v_1^2&v_2^2-u_2v_1^2\\ v_2^2-u_2v_1^2&u_1v_2^2-2u_2v_1v_2\\
  \end{array}
  \right).
\end{equation*}
and (\ref{obs}) is not satisfied: if we denote
$x_1=u_1,\,x_2=u_2,\,x_3=v_1$ and $x_4=v_2$ then the left hand side of
\reff{obs} is given by
%
%
\begin{align*}
  -96x_3(x_4^4-2x_1x_3x_4^3+2x_2x_3^2x_4^2-2x_1x_2x_3^3x_4+x_1^2x_3^2x_4^2+x_2^2x_3^4)
  \,\p^1\wedge\p^2\wedge\p^4
\end{align*}
where the triple wedge is defined by
\begin{equation*}
  \p^i\wedge\p^j\wedge\p^k=\frac16\sum_{\l}\sgn(\l)\,\tt{\l(i)}{\l(j)}{\l(k)}.
\end{equation*}
It follows that in this case the quantized enveloping algebra only
defines a third order deformation quantization.  The choice
$\ph=y^2+xy$ gives another non-zero term; basically any higher order
polynomial leads to an obstruction.  Also the cubic bracket on
$\gl(2)$ (see \cite{lipar}), which is given by
\begin{equation*}
  U=\left(
  \begin{array}{cccc}
    0&x_1^2x_2&x_2x_3^2&x_2x_3(x_1+x_4)\\
    -x_1^2x_2&0&x_2x_3(x_4-x_1)&x_2x_4^2\\
    -x_2x_3^2&x_2x_3(x_1-x_4)&0&x_2x_3^2\\
    -x_2x_3(x_1+x_4)&-x_2x_4^2&-x_2x_3^2&0\\
  \end{array}
  \right)
\end{equation*}
leads to a non-zero obstruction, upon evaluating (\ref{obs}).
Explicitly it is given by
%
%
\begin{multline*}
  96x_2^2x_3(2x_1x_4+x_2x_3)(x_4-x_1)\\
  (x_3\p^1\wedge\p^2\wedge\p^3+(x_4-x_1)\p^1\wedge\p^2\wedge\p^4-
          x_3\p^2\wedge\p^3\wedge\p^4).
\end{multline*}
It follows that for most brackets the enveloping algebra only leads to
a third order deformation.

  \bibliographystyle{amsplain}
\providecommand{\bysame}{\leavevmode\hbox to3em{\hrulefill}\thinspace}

\end{document}